\documentclass{article}

\usepackage{amsmath}
\usepackage{amssymb}
\usepackage{graphicx}
\usepackage{epic,eepic}

\usepackage{color}

\usepackage{amstext,amsthm,amssymb,amsmath}
\usepackage{graphicx}
\usepackage{subfigure}
\usepackage{wrapfig}
\usepackage{fullpage}
\usepackage{color}
\usepackage{multirow}
\usepackage{tabulary}
\usepackage{booktabs}
\usepackage{enumerate}

\usepackage[sort&compress]{}

\newtheorem{remark}{Remark}

\usepackage[utf8]{inputenc}
\usepackage[english]{babel}
\usepackage{hyperref}

\usepackage{nomencl}

\theoremstyle{definition}

\newcommand{\eps}{\varepsilon}

\newcommand{\norm}[1]{\left\|#1\right\|}

\newcommand{\om}{\Omega}

\newcommand\numberthis{\addtocounter{equation}{1}\tag{\theequation}}

\newcommand{\average}[1]{\{\!\!\!\{{#1}\}\!\!\!\}}
\newcommand{\jump}[1]{[\![{#1}]\!]}

\newcommand{\strain}{\mbox{$\boldsymbol{\varepsilon}$}}

\newcommand{\avrg}[2]{{\left \langle #1 \right \rangle}_{#2}}

\hfuzz=\maxdimen
\graphicspath{ {./pics/}}

\title{Adaptive multiscale model reduction with Generalized Multiscale Finite Element Methods}

\author{Eric Chung\thanks{Department of Mathematics, The Chinese University of Hong Kong, Shatin, Hong Kong} \and
Yalchin Efendiev\thanks{Department of Mathematics, Texas A\&M University, College Station, TX 77843} \and Thomas Y. Hou\thanks{Applied and Computational Math. 9-94, California Institute of Technology,  Pasadena 91125}}

\begin{document}
\maketitle

%\tableofcontents
%\newpage
%\section*{Notations}
%\input{notations}
%\input{notations_short}
%\newpage
%\setcounter{page}{1}

\begin{abstract}
In this paper, we discuss a  general multiscale model reduction
framework based on multiscale finite element methods.
We give a brief overview of related multiscale methods. Due to
page limitations, the overview focuses on a  few related
methods and is
not intended to be comprehensive. We present a general adaptive
multiscale model reduction framework,
the Generalized Multiscale Finite Element Method. Besides
the method's basic outline, we discuss some important ingredients
needed for the method's success.
We also discuss several applications. The proposed method allows
performing local model reduction in the presence of high contrast and no scale
separation.

\end{abstract}

\section{Introduction}

\subsection{Objectives}

In this paper, we will present a novel multiscale model reduction technique for solving many challenging problems with multiple scales and high
contrast. Our approach systematically and adaptively adds degrees
of freedom in local regions as needed and goes beyond scale separation cases.  The objectives of the paper are the following:
(1) to demonstrate the main concepts of our unified approach for local multiscale model reduction;
(2) to discuss the method's main ingredients; and
(3) to demonstrate its applications to a variety of challenging multiscale problems.

\subsection{The need for a systematic multiscale model reduction approach}

\begin{figure}[tb]
  \centering
  \includegraphics[width=5in, height=2.5in]{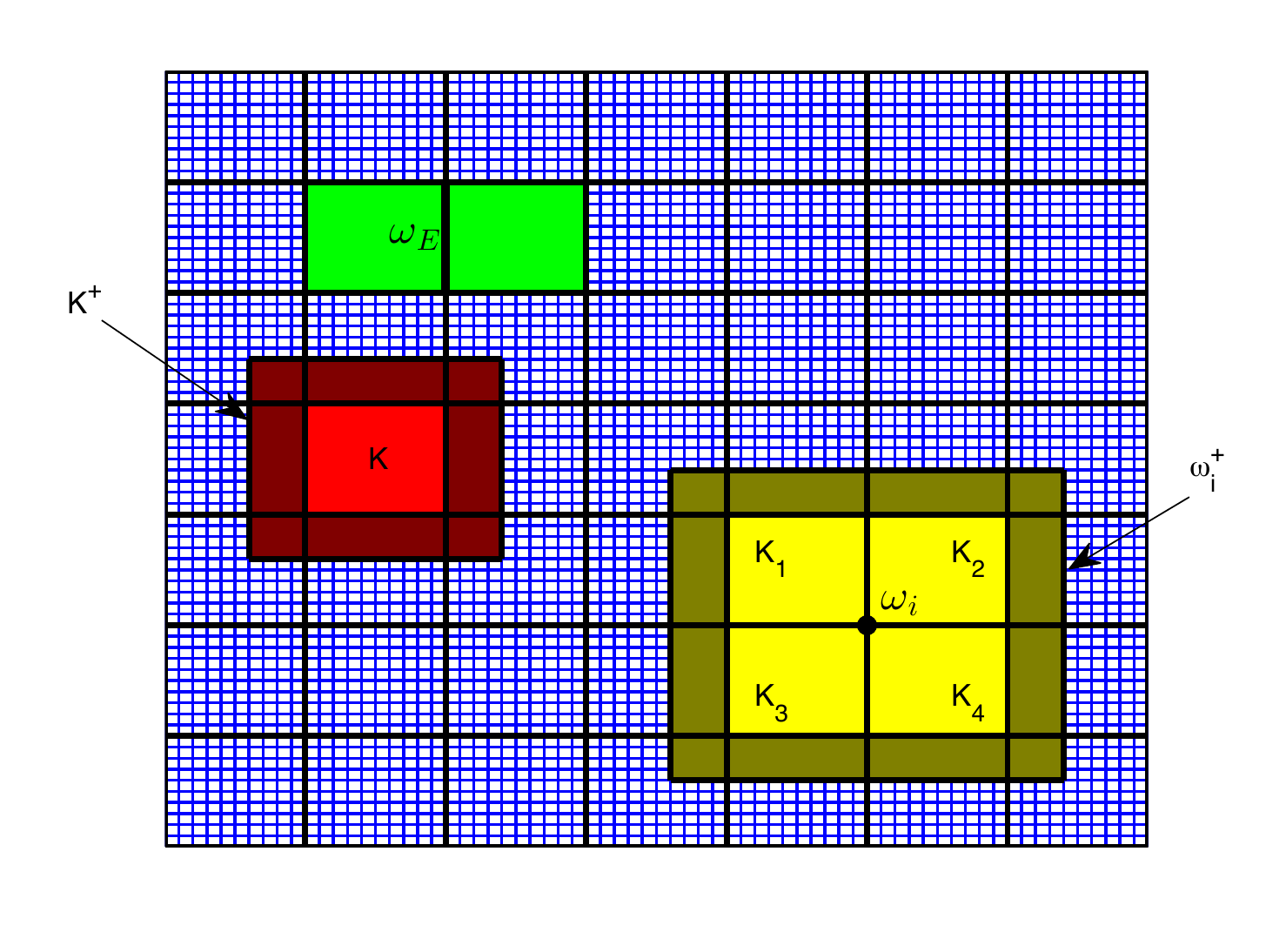}
  \caption{Illustration of fine grid, coarse grid, coarse neighborhood and oversampled domain.}
  \label{schematic_intro}
\end{figure}

{\bf Computational meshes.}
In order to present our approach, we need the notions of coarse and fine meshes,
which are illustrated in Figure \ref{schematic_intro}.
We assume that the computational domain $\Omega$ is partitioned by a coarse grid $\mathcal{T}^H$,
which does not necessarily resolve any multiscale features.
We let $N_c$ be the number of nodes in the coarse grid and $N_e$ be the number of coarse edges.
We use the notation $K$ to represent a generic coarse element in $\mathcal{T}^H$.
To represent multiscale basis functions, we perform
a refinement of $\mathcal{T}^H$ to obtain a fine grid $\mathcal{T}^h$,
with mesh size $h>0$. The fine grid can essentially resolve  all multiscale features of the problem
and we perform the computations of the local basis functions  on the fine grid.

{\bf Scale separation approaches and their limitations.}
Many current approaches for handling complex multiscale problems typically
limit themselves to two or more distinct idealized scales (e.g., pore
and Darcy scales). These include homogenization and numerical homogenization
methods \cite{papanicolau1978asymptotic,  eh09, weh02}.
To demonstrate some of the main concepts, we consider
\begin{equation}
\label{eq:sphase}
\mathcal{L}(u)=f,
\end{equation}
where $\mathcal{L}(u)$ is a differential operator representing the fine-scale
process. One can use many different examples for $\mathcal{L}(u)$. To convey
our main idea, we consider a simple and well-studied heterogeneous diffusion
$\mathcal{L}(u)=-\mbox{div}( \kappa(x)\nabla u)$, where one assumes $\kappa(x)$
to be a multiscale field representing the media properties.
Homogenization and numerical homogenization techniques derive or postulate
macroscopic equations and formulate
local problems for computing the macroscopic parameters.
For example, in the heterogeneous diffusion example, one computes the effective properties
$\kappa^*(x)$ on a coarse grid as a constant tensor, $\kappa_{ij}^*$,
(see Figure \ref{schematic_intro} for illustration
of coarse and fine grids) via
solving local problems (see Section \ref{sec:numhom}).
Using $\kappa^*(x)$, one solves the global problem (\ref{eq:sphase}).
{\it The number of macroscopic parameters represents the effective dimension
of the local solution space}. Consequently, these (numerical homogenization)
approaches cannot
represent many important local features of  the solution space
unless they are identified apriori in a modeling step.
In this paper, we propose a novel approach that avoids
these limitations and determines necessary local degrees of freedom as needed.

{\bf Global model reduction approaches.}
The proposed approaches share some common features with global model
reduction techniques \cite{volkwein05, barrault, ct_POD_DEIM11}, which construct global
basis functions. Our approach uses local
dimension reduction techniques. However, there are many important differences
as we will discuss. First,
global model reduction approaches, though powerful in reducing the degrees of freedom,
lack  local adaptivity and numerical discretization properties
(e.g., conservations of local mass and energy,...)
 that local approaches enjoy. Many successful macroscopic laws (e.g., Darcy's law, and so on)
are possible because the solution space admits a large compression {\it locally}.
We will present several distinct examples to demonstrate this.
For this reason, it is important to construct local multiscale model reduction
techniques that can identify
local degrees of freedom and be consistent with homogenization when there is
scale separation. The proposed
method is a first systematic step in developing such approaches.

{\bf Multiscale Finite Element Methods and some related methods.}
The proposed methods take their origin in Multiscale Finite Element Methods \cite{hw97,jennylt03} and Generalized Finite Element Methods \cite{melenk1996partition}.
 The main idea of MsFEM is to construct local multiscale basis functions,
$\phi^{\omega_i}$ for each coarse block $\omega_i$, and
replace macroscopic equations by using a limited number of basis functions.
More precisely, for each coarse node $i$ (see Figure \ref{schematic_intro}), we
construct multiscale basis function $\phi^{\omega_i}$
 and seek an approximate solution of
(\ref{eq:sphase})
$u_{\text{approx}}=\sum_i c_i \phi^{\omega_i}$.
These approaches motivate our new techniques as MsFEMs are the first methods
that replace macroscopic equation concepts with
carefully designed multiscale basis functions within finite element methods.
Multiscale finite element approaches are shown to be powerful and
have many advantages over numerical
homogenization methods as MsFEMs can recover fine-scale information and
be flexible in terms of gridding. However, these methods do not contain a systematic
way to add degrees of freedom locally and adaptively, which are
 the proposed method's main
contributions.
Approaches, such as variational multiscale methods \cite{hughes98} (see also, \cite{henning2012localized}),
multiscale finite volume \cite{jennylt03}, mortar multiscale method \cite{apwy07},
localized reduced basis \cite{boyaval,ohl12}
can be related to the MsFEM and seek to  approximate  the solution when there is no
scale separation.
%can also be considered as a general multiscale framework. These approaches typically
%correct a coarse grid solution (or stabilize it as it is originally proposed)
%by formulating model equations for the subgrid corrections. These approaches do not contain
%a systematic way of finding reduced solution representation.
Other classes of approaches  built on numerical homogenization methods
(e.g., \cite{ee03,matache2002two,henning2009heterogeneous,fish_book})
are limited to problems with scale separation and when
 macroscopic equations can be formulated.
%lack a framework to add degrees of freedom.
Due to the lack of space, we will briefly discuss a  few of these methods in the paper.

%\begin{figure}
%  \centering
%  \includegraphics[width=0.7 \textwidth]{meshplot}
%  \caption{Schematic of a coarse element and coarse neighborhood}
%  \label{fig:grid}
%\end{figure}

\subsection{The basic concepts of Generalized Multiscale Finite Element Method}

{\bf General idea of GMsFEM.}
In this paper, we will develop novel solution strategies centered
around adaptive multiscale model reduction.
The GMsFEM was first presented in \cite{egh12}
and later investigated in several other papers
(e.g., \cite{galvis2015generalized,Ensemble, eglmsMSDG, eglp13oversampling,calo2014multiscale,chung2014adaptive,chung2015generalizedperforated,chung2015residual,chung2015online}).
The GMsFEM is motivated by our earlier works on designing multiscale
coarse spaces for domain decomposition preconditioners
\cite{ge09_1reduceddim,ge09_2,eglw11}, the studies in MsFEM,
 and the concepts of global
model reduction methods (e.g., \cite{volkwein05,barrault} and references therein), which
employ global snapshots and  Proper Orthogonal Decomposition (POD).
 The proposed
approach is a generalization of the MsFEM and
defines appropriate local snapshots and local spectral decompositions.
The need for such systematic methods will be discussed in the examples of various applications.
Our proposed approaches add local degrees of freedom as needed and provide
{\it numerical} macroscopic
equations for problems without scale separation and high contrast.
Because of the proposed multiscale model reduction's local nature of,
one can adaptively add
the degrees of freedom  based on error estimators and rigorously estimate the
errors.
To our best knowledge, this is one of the
first local systematic multiscale model reduction frameworks that can be easily adopted
for different applications.

{\bf Multiscale basis functions and snapshot spaces.}
The main idea of our multiscale approach is to systematically select
important degrees of freedom for the solution in each coarse
block (see Figure \ref{schematic_intro} for coarse and fine grid illustration).
 More precisely, for each coarse block $\omega_i$ (or $K$),
we identify local multiscale basis functions $\phi_j^{\omega_i}$
($j=1,...,N_{\omega_i}$) and
seek the solution in the span of these basis functions.
 For problems with scale separation, one needs a limited number
of degrees of freedom. However, as the heterogeneities get more
complicated, one needs a systematic approach to find the additional
degrees of freedom.
%The procedure starts with a coarse grid
%layout (see Figure \ref{schematic_intro}).
In each coarse grid, we first build the snapshot space,
$V_{\text{snap}}^{\omega_i}=\text{span}\{\psi_j^{\omega_i}$\}.
 The
choice of the snapshot space depends on the global discretization and
the particular application. The appropriate snapshot space
(1) yields faster convergence, (2) imposes problem relevant
restrictions on the coarse spaces (e.g., divergence free solutions)
and (3) reduces the computational cost of building the offline spaces.
Each snapshot can be constructed, for example, using random boundary
conditions or source terms~\cite{randomized2014}.

{\bf Reducing the degrees of freedom.}
Once we construct the snapshot space $V_{\text{snap}}^{\omega_i}$,
we identify the offline space,
which is a principal component subspace of the snapshot space
and derived based on analysis.
To obtain the offline space, we perform a dimension reduction of the
snapshot space.  This reduction identifies dominant modes, which we use
to construct a multiscale space.  If the snapshot space and the
local spectral problems are appropriately chosen, we obtain reduced
dimensional spaces corresponding to numerical homogenization in the
case of scale separation. When there is no scale separation, we have a
constructive method to add the necessary extra degrees of freedom that
capture the relevant interactions between scales.

{\bf Adaptivity and nonlinearities.}
The algorithmic framework proposed above is general and can be used
for different multiscale, high contrast, and perforated problems~\cite{
  eric-2012, chung2015generalizedperforated, galvis2015generalized,
  calo2014multiscale}.  The multiscale basis functions are constructed
locally and using an adaptive criteria.  Thus, in different regions,
we expect a different number of basis functions depending on the local
features of the problem, such as heterogeneities and high contrast.
For example, in the regions with scale separation, we expect only a
limited number of degrees of freedom.  The adaptivity will be
achieved using error indicators.  In nonlinear problems, one performs
local nonlinear interpolation to approximate the Jacobian or other
nonlinear terms.  The above multiscale procedure can be complemented
with online basis functions which helps to converge to the fine-scale
solution by constructing multiscale basis functions in the
simulations~\cite{ chung2015residual, chung2015online} (also, see
Figure~\ref{fig:mmr1} for the ingredients of proposed multiscale
method).

{\bf Limitations.}
Though the proposed approaches can be used in many applications,
it is limited to problems where the solution space locally has a low
dimensional structure.
If the latter is not the case, our approaches will use many degrees
of freedom and the computational gain may not be significant in those regions.
Many multiscale application problems we have encountered can benefit
from local adaptive multiscale model reduction.

\subsubsection{Methodological ingredients of GMsFEM}

\begin{figure}[t]
\centering
\includegraphics[width=4in, height =3in]{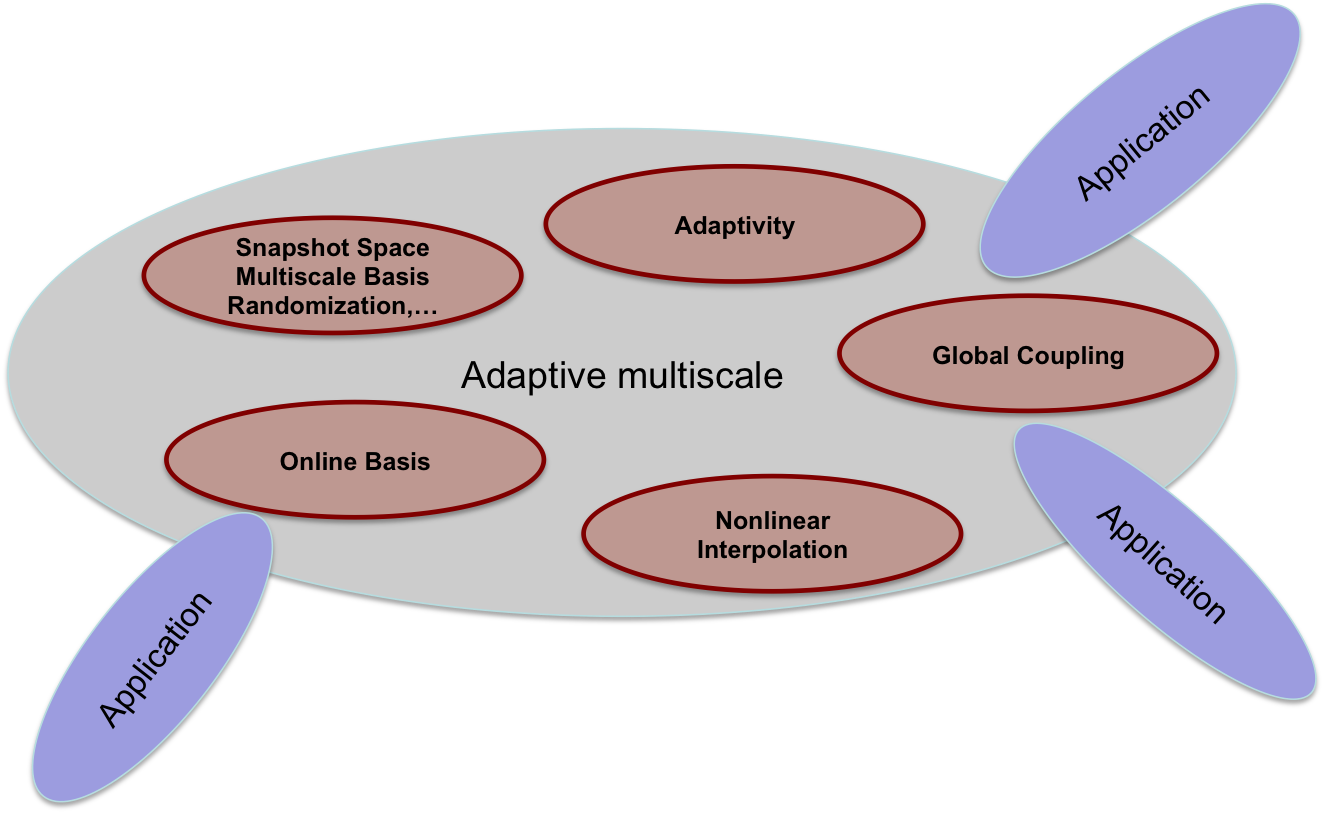}
\caption{Multiscale model reduction. Ingredients and Applications.}
\label{fig:mmr1}
\end{figure}

Our computational framework will rely on several
important ingredients (see Figure~\ref{fig:mmr1})
that are general for various discretizations (such as mixed methods,
discontinuous Galerkin, and etc,) and applications. These ingredients
include (1)
a procedure for identifying local snapshot spaces and multiscale basis,
(2) developing global coupling mechanisms for multiscale basis
functions, (3) adaptivity strategies, (4) local nonlinear interpolation
tools, and (5) online basis functions.
These ingredients are building blocks that are needed for
constructing an accurate and robust multiscale model reduction.
We will motivate the use of
multiscale model reduction by showing how it can take us beyond
conventional macroscopic modeling and discuss some relevant
ingredients.

%\subsection{A promising framework}

The proposed framework provides a promising tool and capability
for solving a large class of multiscale
problems without scale separation and high contrast.
This framework is tested and
applied in various applications where we have followed the general
framework to construct local multiscale spaces.
We identify main ingredients of the framework, which can be further
investigated for speed-up and accuracy.

\subsection{Applications and organization of the paper}

We apply the proposed adaptive multiscale framework to several
challenging and distinct
applications where one deals with a rich hierarchy of
scales. These include (1) flows in heterogeneous porous media (2) diffusion
in fractured media (3) multiscale processes in perforated regions (4) wave
propagation in heterogeneous media (5) nonlinear
diffusion equations in heterogeneous media.
We will discuss some of the applications here.

%\subsection{Organization of the paper}

In Section \ref{sec:numhom}, we describe numerical homogenization
and its limitations.
In Section \ref{sec:msfem}, we discuss Multiscale Finite Element Method.
In Section \ref{sec:gmsfem_basics}, we describe the GMsFEM and its basic
principles.
In Section \ref{sec:adaptivity}, adaptive strategies for the GMsFEM
are described. Section \ref{sec:online} is devoted to the construction
of the online basis functions.  We describe the selected global
couplings in Section \ref{sec:selected_discretization}.
The GMsFEM using sparsity in the snapshot space is discussed
in Section \ref{sec:sparsity}. We discuss the GMsFEM for space-time heterogeneous
problems in Section \ref{sec:spacetime}.
In Section \ref{sec:perforated}, we discuss the GMsFEM for problems in
perforated domains.
 We present some selected applications
in Section \ref{sec:applications}.
In Section \ref{sec:summary}, we discuss some remaining
aspects of the GMsFEM, such as applications to nonlinear
problems and its use in global model reduction techniques.

\section{A brief introduction to numerical homogenization}
\label{sec:numhom}

\subsection{Numerical homogenization}

The {\it main idea of numerical homogenization} is to identify
the homogenized coefficients in each coarse-grid block.
The basic underlying principle is to compute these
upscaled quantities such that they preserve some averages
for a given set of local boundary conditions.
We discuss it on an example.

We consider a simple example
\begin{equation}
\frac{\partial}{\partial x_i} \left ( \kappa_{ij}(x) \frac {\partial}{\partial x_j} u \right ) = f,
\end{equation}
where $\kappa_{ij}(x)$ is a heterogeneous field.
The objective is to define an upscaled (or homogenized)
conductivity for each coarse block without directly using the
periodicity. Following the homogenization technique, the local problems
are solved in each coarse block
\begin{equation}
\label{eq:numhomeq}
\frac{\partial}{\partial x_i} \left ( \kappa_{ij}(x) \frac {\partial}{\partial x_j} \mathcal{N}_l\right ) = 0 \ \text{in} \ K,
\end{equation}
where $K$ is a coarse block.
The boundary condition needs to be chosen such that it represents
the local heterogeneities. We limit ourselves to
Dirichlet boundary conditions
${\mathcal{N}}_l=x_l$ on $\partial K$.
One can use other boundary conditions \cite{weh02}.

We note that if the local problem is homogenized and
the homogenized coefficients are constant, then the solution
of the homogenized equation is
$\mathcal{N}_l^*=x_l$ in $K$.
The upscaled coefficients $\kappa^{*,nh}_{ij}$ can be defined by averaging the fluxes:
\begin{equation}
\int_K \kappa_{ij}^{*,nh} \frac {\partial}{\partial x_j} \mathcal{N}_l^* =
\int_K \kappa_{ij}(x) \frac {\partial}{\partial x_j} \mathcal{N}_l.
\end{equation}
The motivation behind this upscaling is to state that the average flux
response for the fine-scale local problem with prescribed boundary
conditions is the same as that for the upscaled solution.
If we take $\mathcal{N}_l^*=x_l \ \text{in}\ \ K$, we have
\begin{equation}
\kappa_{il}^{*,nh}  =
{1\over |K|} \int_K \kappa_{ij}(x) \frac {\partial}{\partial x_j} \mathcal{N}_l.
\end{equation}

%\begin{figure}[tbp]
%\centering
%\includegraphics[width=4in, height=2.in]{numhom.eps}
%\caption{Illustration of numerical homogenization}
%\label{fig:numhom}
%\end{figure}

\begin{remark}[Convergence]
Assuming that
$\kappa_{ij}(x)=\kappa_{ij}(x/\epsilon)$,
one can study the convergence of the numerical homogenization technique and
 show (\cite{weh02}) that
\[
|\kappa_{il}^{*,nh}- \kappa_{il}^*|\leq C \sqrt {\epsilon\over H},
\]
where $\kappa_{il}^{*,nh}$ is
the numerical homogenization and  $\kappa_{il}^*$ is the
``correct''
 homogenized
coefficients.

\end{remark}

\begin{remark}[Oversampling]
The convergence analysis suggests that the boundary layers due to
artificial linear boundary conditions, $\mathcal{N}_l=x_l$,
cause the resonance errors. To
reduce these resonance errors,
oversampling technique has been proposed. The main
idea of this method is to solve local problems in a larger domain.
In oversampling method \cite{hw97}, the local problem that is analogous to
(\ref{eq:numhomeq}) is solved in a larger domain (see Figure
\ref{schematic_intro}). In particular, if we denote the large
domain by $K^+$ while the target coarse block by $K$ (see Figure
\ref{schematic_intro}), then
\begin{equation}
\label{eq:numhomeqovs}
\frac{\partial}{\partial x_i} \left ( \kappa_{ij}(x) \frac {\partial}{\partial x_j} \mathcal{N}_l^{ovs}\right ) = 0 \ \text{in} \ K^+.
\end{equation}
For simplicity, we can use
$\mathcal{N}_l^{ovs} = x_l \ \text{on}\ \partial K^+$.
The upscaled conductivity $\kappa^{*,ovs}_{ij}$ is computed by equating the fluxes on
the target coarse block:
\begin{equation}
\int_K \kappa_{ij}^{*,ovs} \frac {\partial}{\partial x_j} \mathcal{N}_l^{*,ovs} =
\int_K \kappa_{ij}(x) \frac {\partial}{\partial x_j} \mathcal{N}_l^{ovs}.
\end{equation}
Noting that the domain $K^+$ is only slightly larger, we can still claim
that
$ \mathcal{N}_l^{*,ovs}= x_l$,
and thus,
\begin{equation}
\kappa_{il}^{*,ovs} =
{1\over |K|} \int_K \kappa_{ij}(x) \frac {\partial}{\partial x_j} \mathcal{N}_l^{ovs}.
\end{equation}
The advantage of this approach is that one reduces the effects of the
oscillatory boundary conditions. One can show (\cite{weh02})
that the error for the residual
is small and scales as $\epsilon/H$ instead of $\sqrt{\epsilon/H}$.

\end{remark}

%\begin{figure}[tbp]
%\centering
%\includegraphics[width=4in, height=2.5in]{oversampling.eps}
%\caption{Coarse and fine grid}
%\label{fig:oversampling}
%\end{figure}

\subsection{Limitations}

The above numerical homogenization concepts (for linear and nonlinear
problems) are often used to derive
macroscopic equations.  These techniques assume that the solution
space has a reduced dimensional structure.  For example, as we can
observe from the above example and derivations, the solution in each
coarse block is approximated by three local fields (in 3D).  In fact,
this ``effective'' dimension is related to the number of elements
(constants) in the effective properties.  In general, for many
macroscopic equations, one implicitly assumes the local effective
dimension for the solution space related to the number of
macroscopic parameters.  In many cases, the assumptions on the limited
effective dimension of microscale problem break down and the local
solution space may need more degrees of freedom.  This requires
general approaches, where we do not rely on macroscale equations, and
approximate the solution space via multiscale basis functions.  Next,
we give a brief overview of these concepts.

\section{Multiscale Finite Element Method}
\label{sec:msfem}

In this section, we will give a brief overview of MsFEM as a method
for solving a problem on a coarse grid.
MsFEMs consist of two major ingredients: (1) multiscale basis functions
and (2) a  global numerical formulation
which couples these multiscale basis functions.
Multiscale basis functions are designed to capture the
fine-scale features of the solution.
Important multiscale features of the solution need to be incorporated
into these localized basis functions which contain
information about the scales which are smaller (as well as
larger) than the local numerical
scale defined by the basis functions.
In particular, we need to incorporate the features of the solution
that can be localized and use additional basis functions to
capture the information about the features that need to be separately included
in the coarse space.
A global
formulation couples these basis functions to provide an accurate
approximation of the solution.

As before, we consider the second order elliptic equations with heterogeneous
coefficients
\begin{equation}
\label{main:red}
\mathcal{L}(u)=f,
\end{equation}
where $\mathcal{L}(u)=-\frac{\partial}{\partial x_i} \left ( \kappa_{ij}(x) \frac {\partial}{\partial x_j} u \right )$ with appropriate boundary conditions.
We seek multiscale basis functions supported in each domain $\omega_j$, denoted
by $\phi^{\omega_j}$. Then, the coarse-grid solution is represented by
\[
u_H=\sum_i c_i \phi^{\omega_i},
\]
where $c_i$ are determined from
\begin{equation}
\label{coarse2}
a(u_H,v_H)=(f,v_H),\quad \text{for all}\  v_H\in V_0.
\end{equation}
In the above formulation, we define
 $V_0=\text{span}\{ \phi^{\omega_i} \}$,
\begin{equation}\label{eq:def:a}
a(u,v)=\int_\Omega
\kappa(x)\nabla u(x)\nabla v(x) \quad \text{and} \quad
(f,v)=\int_\Omega f(x)v(x).
\end{equation}

One can also view MsFEM in the discrete setting.
Assume that the basis functions are defined on a fine
grid as $\Phi^{\omega_i}$ with $i$ varying from
$1$ to $N_c$, where $N_c$ is the number of multiscale basis functions.
Given   coarse-scale basis functions,
the coarse space is given by
\begin{equation}\label{eq:def:V0}
V_0=\mbox{span} \{ \Phi^{\omega_i}\}_{i=1}^{N_c},
\end{equation}
and the coarse matrix is given by  $A_H=R_0^TAR_0$ where $A$ is the fine-scale stiffness matrix and
\[
R_0=[\Phi^{\omega_1},\dots,\Phi^{\omega_{N_c}}].
\]
Here $\Phi^{\omega_i}$'s are discrete coarse-scale basis functions
defined on a fine grid
(i.e., column vectors).
 Multiscale finite element solution is the finite
element projection of the fine-scale solution into the space
$V_0$. More precisely, multiscale solution $U_H$
is given by
\[
A_H u_H=f_H,
\]
where $f_H=R_0^T b$. Next, we discuss some coarse spaces.

%\begin{figure}[tbp]
%\centering
%\includegraphics[width=3in, height=2in]{basis.eps}
%\caption{Illustration of some multiscale basis functions.}
%\label{ill_basis_many}
%\end{figure}

%\begin{figure}[tbp]
%\centering
%\includegraphics[width=3in, height=2in]{ill_omega.eps}
%\caption{Schematic description of coarse regions.}
%\label{ill_omega}
%\end{figure}

%We discuss some coarse spaces
%constructed to capture
%the fine-scale features of the solution. In particular, we will
%{\it first consider coarse spaces where there is only one function
%per coarse node} $y_i$ is defined.
%For this reason, we will use the notation $\chi_i$.
%Further, we discuss how
%these basis functions can be complemented
%so that our coarse-scale approximation converges
%to the fine-scale solution rapidly.

{\bf Linear boundary conditions.}
%\label{sec:lin_boundary_conditions}
Let  $\chi_i^0$ be
the standard piecewise linear or piecewise polynomial basis function supported in
$\omega_i$.
We define multiscale finite element basis functions that
coincide with  $\chi_i^0$ on the boundaries of the coarse partition.
In particular,
\begin{eqnarray}
\mbox{div}(\kappa\nabla\chi_i^{ms})=0\ \ \mbox{in }K\in \omega_i,\quad
\chi_i^{ms}=\chi_i^0\ \ \mbox{in }\partial K,\ \ \forall\  K \subset \omega_i,
\label{e}
\end{eqnarray}
where $K$ is a coarse grid block within $\omega_i$.
Note that multiscale basis functions coincide with standard finite element
basis functions on the boundaries of coarse grid blocks, while are oscillatory
in the interior of each coarse grid block.

\begin{remark}
We would like to remark that the
MsFEM formulation allows one to take advantage of
scale separation.
In particular, $K$ can be chosen
to be a volume smaller than the coarse grid
and the integrals in the stiffness matrix computations need to be re-scaled
(see \cite{eh09} for discussions).

\end{remark}

\begin{remark}[The relation between MsFEM and numerical homogenization]
%The relation between the MsFEM and numerical homogenization can be investigated
%(\cite{eh09}).
It can be shown (\cite{eh09}) that the MsFEM with one basis function
on triangular elements yields the same coarse-grid stiffness matrix
as the numerical homogenization.This is due to the fact that the numerical
homogenization uses the local solutions of PDEs as in the MsFEM. However,
the MsFEM has several advantages, which include  fine-scale information
recovery, adaptivity based on the residual,
the use of global information and so on (\cite{eh09}).
\end{remark}

\begin{remark}[The relation between MsFEM and some other multiscale techniques]
The relation between the MsFEM and other multiscale methods is discussed
in \cite{eh09}. It can be shown that the MsFEM can use an approximation
of multiscale basis functions if there is a periodicity. As a result,
the MsFEM can recover a similar approximation and the computational cost
as Heterogeneous Multiscale Method (\cite{ee03}) for elliptic equations
 in the presence of scale
separation. The relation of the
MsFEM and variational multiscale method (\cite{hughes98}) is also discussed
in \cite{eh09}. The variational multiscale method (VMM) recovers the fine-scale
information via the equation for the residual. The localization of this equation
and the choice of the initial multiscale spaces are important for the VMM.
It can be shown
that with a simple choice, one can recover the MsFEM.
On the other hand, one can use more sophisticated coarse spaces
and the residual recovery, to improve the accuracy of the VMM \cite{calo2011note}.
\end{remark}

{\bf Oversampling technique.}
Because of linear boundary conditions, the basis functions do not
capture
the fine-scale features of the solution along the boundaries.
This can lead to large errors. When coefficients have a single
physical scale $\epsilon$, it has been
shown that the error (see \cite{hwc99}) is
proportional to $\epsilon/H$, and thus can be large
when $H$ is close to $\epsilon$.
Motivated by such examples, Hou and Wu in \cite{hw97}
proposed an oversampling technique for multiscale finite element
method. Specifically, let $\psi^{+,\omega_i}$
be the basis functions satisfying
the homogeneous elliptic equation in the larger
domain $K^{+} \supset K$ (see Figure \ref{schematic_intro}).
We then form the actual basis $\phi^{+,\omega_i}$ by linear combination of
$\psi^{+,\omega_i}$,
$\phi^{+,\omega_i} = \sum_{j} \alpha_{ij} \psi^{+,\omega_j}$, and
restricting them to $K$.
The coefficients $\alpha_{ij} $ are determined by
condition $\phi^{+,\omega_i} (x_j ) = \delta_{ij} $, where $x_j$ are nodal points.
Other conditions can also be imposed (e.g., $\alpha_{ij}$ are determined
based on homogenized parts of $\psi^{+,\omega_i}$).
Note that this method is non-conforming.
One can also multiply the oversampling  functions
by linear basis functions to restrict them onto $\omega_i$ and
have a conforming method.
Numerical results and more discussions on oversampling methods
can be found in \cite{eh09}.
Many other boundary conditions are introduced and analyzed
in the literature. For example, reduced boundary conditions
are found to be efficient in many porous media applications
(\cite{jennylt03}).

{\bf The use of limited global information.}
Previously, we discussed multiscale methods which employ local information in
computing basis functions with the exception of
energy minimizing basis functions.  The accuracy of these approaches
depends on local boundary conditions. Though effective in many cases,
multiscale methods that only use local information may not accurately
capture the local features of the solution. In a number of
previous papers, multiscale methods that employ limited
global information is introduced.
The main idea of these multiscale methods is to incorporate some
fine-scale information about the solution that can not be
computed locally and that is given globally.
More precisely, in these approaches,
 we assume that the solution can be represented
by a number of fields  $p_1$,..., $p_N$, such that
\begin{equation}
\label{assump12}
u \approx G(p_1,...,p_N),
\end{equation}
where $G$ is sufficiently smooth
function, and $p_1$,.., $p_N$ are global fields.
These fields  typically
contain the essential information about the heterogeneities
at different scales and can also be local fields defined in larger domains.
In the above assumption (\ref{assump12}), $p_i$ are solutions of elliptic equations. These global fields are used to construct multiscale basis
functions (often multiple basis functions per a coarse node).
Finding $p_1$, ..., $p_N$, in general, can be a difficult task
and we refer to Owhadi and Zhang \cite{oz07} as well as
to \cite{eghe05} where various choices of global information are proposed.

\section{Generalized Multiscale Finite Element Method. Basic Concepts}
\label{sec:gmsfem_basics}

\subsection{Overview}

Previous multiscale finite element concepts focus on constructing
one basis function per node. As we discussed earlier
 these approaches
are similar to numerical homogenization (upscaling), i.e.,
one can use one basis function to localize the effects of
local heterogeneities. However, for many complex heterogeneities,
multiple basis functions are needed to represent the local solution
space. For
example, if the coarse region contains several high-conductivity
regions, one needs multiple multiscale basis functions to represent the local
solution space. In this section,
we discuss how to systematically construct multiscale basis functions
in a general framework
Generalized Multiscale Finite Element
Methods.

\begin{figure}[tbp]
\centering
\includegraphics[width=5in, height=2.5in]{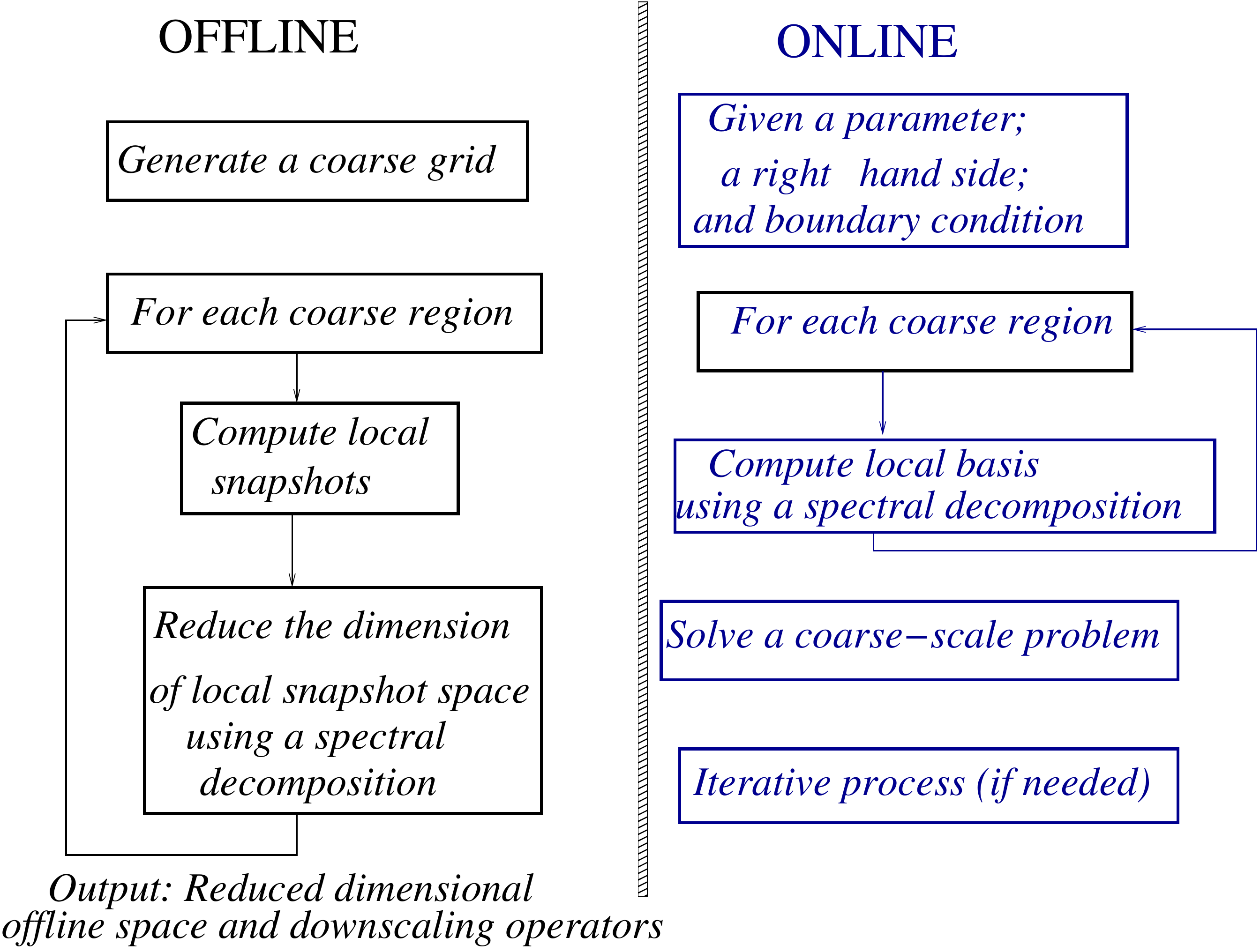}
\caption{Flow chart}
\label{scheme}
\end{figure}

Generalized Multiscale
Finite Element Method (GMsFEM)
incorporates complex input space information
and the input-output relation.
 It systematically enriches the
coarse space through our local construction. Our approach, as in many
multiscale and model reduction techniques, divides the computation into
two stages (see Figure \ref{scheme}): offline and online. In the offline stage, we
construct a small dimensional space that can be efficiently
used in the online stage to construct multiscale basis functions.
These multiscale basis functions can be re-used for any input parameter
to solve the problem on a coarse-grid. Thus, this provides a substantial
computational saving in the online stage. Below, we present an outline
of the algorithm and will discuss its main ideas on the example
of (\ref{main:red}).

{\bf Offline computations:}
 (1) Coarse grid generation;
(2) Construction of snapshot space that will be used to
compute an offline space.
(3) Construction of a small dimensional offline space
by performing dimension reduction in the space of snapshots.

{\bf Online computations:}
(1) For each input parameter, compute multiscale basis functions
(for parameter dependent problems);
(2) Solution of a coarse-grid problem for any force
term and boundary condition;
(3) Iterative solvers, if needed.

%\begin{figure}[tbp]
%\centering
%\includegraphics[width=4in, height=2.in]{chart_IO.eps}
%\caption{Flow chart}
%\label{scheme1}
%\end{figure}

\subsection{Examples of snapshot spaces. Oversampling and non-oversampling}

The snapshot space, denoted by $V_{H,\text{snap}}^{\omega_i}$ for a generic domain
$\omega_i$,
 is defined for each coarse region
$\omega_i$
 (with elements of the space denoted $\psi_l^{\omega_i}$) and can approximate
 the solution space with a
prescribed accuracy.
 The appropriate snapshot space
(1) yields faster convergence, (2) imposes problem relevant
restrictions on the coarse spaces (e.g., divergence free solutions)
and (3) reduces the computational cost of building the offline spaces.
In particular, we emphasize that the use of oversampling in the snapshot
spaces can improve the convergence.

\subsubsection{All fine-grid functions}

We can use local fine-scale spaces consisting of fine-grid basis
functions within a coarse region. More precisely,
the snapshot vectors consists of unit vectors defined on
a fine grid within a coarse region.
 In this case, the offline spaces
will be computed on a fine grid directly (\cite{ge09_2}).
%The local fine-grid space has an advantage  if the dimension of the local
%fine spaces  is comparable to the dimension of
%$V_{H,\text{snap}}^{\omega_i}$.
% computed by solving local problems
%as in the first choice.

\subsubsection{Harmonic extensions}

This choice of snapshot space consists of harmonic extension
of fine-grid functions defined on the boundary of $\omega_i$.
More precisely, for each fine-grid function, $\delta_l^h(x)$,
which is defined by
$\delta_l^h(x_k)=\delta_{l,k},\,\forall  x_k\in \textsl{J}_{h}(\omega_i)$, where
$\textsl{J}_{h}(\omega_i)$ denotes the set of fine-grid boundary nodes on $\partial\omega_i$, we obtain a snapshot function $\psi_l^{\omega_i}$ by
\[
\mathcal{L}( \psi_{l}^{\omega_i})=0\ \ \text{in} \ \omega_i
\]
subject to the boundary condition, $ \psi_{l}^{\omega_i}=\delta_l^h(x)$
on $\partial\omega_i$.
Here $\delta_{l,k} = 1$ if $l=k$ and $\delta_{l,k}=0$ if $l\ne k$.

\subsubsection{Oversampling approaches}

To describe the oversampling approach, we consider the harmonic
extension as described above in the oversampled region
(see, e.g., Figure \ref{schematic_intro}).
This choice of snapshot space consists of harmonic extension
of fine-grid functions defined on the boundary of $\omega_i^{+}$.
More precisely, for each fine-grid function, $\delta_l^h(x)$,
which is defined by
$\delta_l^h(x_k)=\delta_{l,k},\,\forall x_k\in \textsl{J}_{h}(\omega_i^{+})$, where $\textsl{J}_{h}(\omega_i^{+})$ denotes the set of fine-grid boundary nodes on $\partial\omega_i^{+}$,
we obtain a snapshot function $\psi_l^{+,\omega_i}$ by
\[
\mathcal{L}( \psi_{l}^{+, \omega_i})=0\ \ \text{in} \ \omega_i^{+}
\]
subject to the boundary condition, $ \psi_{l}^{+, \omega_i}=\delta_l^h(x)$
on $\partial\omega_i^{+}$.

\subsubsection{Randomized boundary conditions}

In the above construction of snapshot vectors, many local problems
are solved. This is not necessary and one can only solve a relatively small
number of snapshot vectors. The number of snapshot vectors
is defined by the number of offline multiscale basis functions.
In this case, we will use random boundary conditions.
%We describe this procedure in more details later on.

More precisely, for each fine-grid function, $r_l^h(x)$,
which is defined by
$r_l^h(x_k)=r_{l,k},\,\forall x_k\in \textsl{J}_{h}(\omega_i^{+})$,
where $r_{l,k}$ are random numbers,
we obtain a snapshot function $\psi_l^{+,\omega_i}$ by
\[
\mathcal{L}( \psi_{l}^{+, \omega_i})=0\ \ \text{in} \ \omega_i^{+}
\]
subject to the boundary condition, $ \psi_{l}^{+, \text{snap}}=r_l^h(x)$
on $\partial\omega_i^{+}$.

\subsection{Offline spaces}

The offline space, denoted by $V_{H,\text{off}}^{\omega_i}$ for a generic domain
$\omega_i$,  is defined for each coarse region
$\omega_i$
 (with elements of the space denoted $\phi_l^{\omega_i}$) and
used to approximate the solution.
{  The offline space is constructed by performing a spectral decomposition
in the snapshot space and selecting the dominant eigenvectors
(corresponding to the smallest eigenvalues).}
{ The choice of the spectral problem is important for
the convergence and is derived from the analysis as it is described below.
The convergence rate of the method is proportional to $1/\Lambda_*$, where
$\Lambda_*$ is the smallest eigenvalue among all coarse blocks whose
corresponding
eigenvector is not included in the offline space.
Our goal is to select the local spectral problem
to remove as many small eigenvalues as possible so that we can obtain smaller
%there are fewer very
%small eigenvalues and, consequently, with smaller
dimensional coarse spaces to achieve a higher accuracy.}

\subsubsection{General concept and example}

 The construction of the offline space requires solving an
appropriate local spectral
eigenvalue problem.
The local spectral problem is derived from the analysis.
A first step in the analysis is to decompose the energy functional
corresponding to the error into coarse subdomains.
For simplicity, we denote the energy functional corresponding to
the domain $\Omega$ by $a_\Omega(u,u)$, e.g.,
$a_\Omega(u,u) = \int_\Omega \kappa \nabla u\cdot \nabla u$.
Then,
\begin{equation}
\begin{split}
a_\Omega(u-u_H,u-u_H)\preceq
\sum_\omega a_\omega(u^\omega-u_H^\omega,u^\omega-u_H^\omega),
\end{split}
\end{equation}
where $\omega$ are coarse regions ($\omega_i$),
$u^\omega$ is the localization of the solution using
the snapshot vectors defined for $\omega$
and $u_H^\omega$ is the component of the solution $u_H$ spanned by the basis in $\omega$.

The local spectral decomposition is chosen to bound
$ a_\omega(u^\omega-u_H^\omega,u^\omega-u_H^\omega)$.
We seek
the subspace $V_{H,\text{off}}^{\omega}$
(since $u^\omega \in V_{H,\text{snap}}^{\omega}$)
such that for any
$\psi\in V_{H,\text{snap}}^{\omega}$,
there exists
$\psi_0\in V_{H,\text{off}}^{\omega}$ with,
\begin{equation}
\label{eq:off1}
a_{\omega}(\psi-\psi_0,\psi-\psi_0)\preceq {\delta}
s_{\omega}(\psi-\psi_0,\psi-\psi_0),
\end{equation}
where
$s_{\omega}(\cdot,\cdot)$ is an auxiliary bilinear form,
which needs to be chosen and $\delta$ is a threshold related to eigenvalues.
 First, we would like the fine-scale
solution to be bounded
in $s_{\omega}(\cdot,\cdot)$.
 We note that, in computations, (\ref{eq:off1})
involves solving a generalized eigenvalue problem
with a mass matrix defined using $s_{\omega}(\cdot,\cdot)$ and the basis functions are selected based on
dominant eigenvalues as described above.  The threshold value $\delta$ is
chosen based on the eigenvalue distribution.
Secondly, we would like
the eigenvalue problem to have a fast decay in the spectrum
(this typically requires using oversampling techniques).
Thirdly, we would like to
bound
\[
\sum_\omega s_\omega(u^\omega-u_H^\omega,u^\omega-u_H^\omega)\preceq \widetilde{a}_\Omega(u,u),
\]
where $\widetilde{a}_\Omega(u,u)$ can be bounded independent of physical
parameters and the mesh sizes. In this step, one can need energy minimizing
snapshots (See Remark \ref{rem:energy_min} and \cite{chung2015_EMO}).

Next, we discuss a choice for $a_\omega(\cdot,\cdot)$ and
$s_\omega(\cdot,\cdot)$.
Recall that for (\ref{eq:off1}), we need a
a local spectral problem, which is to find a real number $\lambda$ and $v \in V_{H,\text{snap}}^{\omega_i}$ such that
\begin{equation}\label{spectralProblem_GMsFEM}
a_{\omega_i}(v, w) = \lambda s_{\omega_i}(v, w), \qquad \forall w \in V_{H,\text{snap}}^{\omega_i},
\end{equation}
where $a_{\omega_i}$ is a symmetric non-negative definite bilinear operator
and $s_{\omega_i}$ is a symmetric positive definite bilinear operators defined on $V_{H,\text{snap}}^{\omega_i} \times V_{H,\text{snap}}^{\omega_i}$.
Based on our analysis, we can choose
\begin{equation*}
\begin{aligned}
%a_{\omega_i}(v, w) = \int_{\omega_i} \kappa \nabla (\chi_i^{ms} v) \cdot \nabla (\chi_i^{ms} w), \
a_{\omega_i}(v, w) = \int_{\omega_i} \kappa \nabla v \cdot \nabla w, \
s_{\omega_i}(v, w) = \int_{\omega_{i}} \widetilde{\kappa}  v   w,
\end{aligned}
\end{equation*}
where $\widetilde{\kappa} =  \sum_{i=1}^{N_c}  \kappa\nabla \chi_i^{ms}  \cdot \nabla \chi_i^{ms}$
and $\chi_i^{ms}$ are multiscale basis functions
(see (\ref{e})).
We let $\lambda_j^{\omega_i}$ be the eigenvalues
of (\ref{spectralProblem_GMsFEM}) arranged in ascending order.
We will use the first $l_i$ eigenfunctions to construct our
offline space $V_{H,\text{off}}^{\omega_i}$.
{ The choice of the eigenvalue problem
is motivated by the convergence analysis.
 The convergence rate is proportional to $1/\Lambda_*$ (with
the proportionality constant depending on $H$), where $\Lambda_*$ is
the smallest eigenvalue that
the corresponding eigenvector is not included in the offline space.
We would like to
remove as many small eigenvalues as possible.
The eigenvectors of the corresponding small eigenvalues
represent important features of the solution space.
% have fewer very small
% eigenvalues, for which the corresponding eigenvectors
%represent important features of the solution space.
By choosing the multiscale basis functions, $\chi_i^{ms}$, in the construction
of local spectral problem, some localizable important features are taken
into consideration through $\chi_i^{ms}$
and, as a result, we have fewer small eigenvalues
(see \cite{ge09_2} for more discussions).}

For constructing conforming multiscale basis functions,
the selected eigenfunctions are multiplied
 by the partition of unity
functions, such as $\chi_i^{ms}$ or $\chi_i^0$
(cf., \cite{melenk1996partition}).
The multiplication by the partition of unity functions modifies
 the multiscale
nature of the selected multiscale eigenfunctions.
We will discuss some other important
discretizations such as mixed, discontinuous Galerkin discretizations, which
can be more suitable for the coupling of the multiscale eigenfunctions and
which avoid the multiplication by the partition of
unity functions.
The global offline space
$V_{H,\text{off}}$ is formed as the union of all
$V_{H,\text{off}}^{\omega_i}$.
Once the offline space is constructed,
we solve (\ref{coarse2}) and find
$u_H=\sum c_{i,j} \phi_j^{\omega_i}\in V_{H,\text{off}}$ such that
\begin{equation}
\label{eq:GMsFEM_global}
a(u_H,v_H)=(f,v_H), \ \forall v_H\in V_{H,\text{off}}.
\end{equation}

\subsubsection{An implementation view }

Next, we present an implementation view of the local spectral decomposition.
The space formed by the snapshot vectors is
\[
V_{H,\text{snap}}^{\omega_i}
= \text{span}\{ \Psi_{l}^{\omega_i}:~~~ 1\leq l \leq L_i \},\ \
V_{H,\text{snap}}^{+,\omega_i}=
\text{span}\{ \Psi_{l}^{+, \omega_i}:~~~1\leq l \leq L_i^+ \}
\]
for each  coarse neighborhood $\omega_i$
and
for each oversampled coarse neighborhood $\omega_i^+$, respectively, where $L_i$ and $L_i^+$
are the dimensions of the snapshot spaces $V_{H,\text{snap}}^{\omega_i}$ and $V_{H,\text{snap}}^{+,\omega_i}$.
We note that in the case when $\omega_i$ is adjacent to
the global boundary, no oversampled domain is used.
We can put all snapshot functions using a matrix representation
$$
R_{\text{snap}}^{\omega_i} = \left[ \Psi_{1}^{\omega_i}, \ldots, \Psi_{L_i}^{\omega_i} \right],\quad
R_{\text{snap}}^{+,\omega_i} = \left[ \Psi_{1}^{+,\omega_i}, \ldots, \Psi_{L_i^+}^{+,\omega_i} \right],
$$
%where $\Psi_j^{\text{snap}}$ denotes the restriction of $\Psi_j^{\text{+,snap}}$ to $\omega_i$, and
%where
%$M_{\text{snap}}$ denotes the total number of functions to keep in the snapshot matrix construction.

The local spectral problems (\ref{spectralProblem_GMsFEM}) can be written in a matrix form as
\begin{equation}
A^{\omega_i} \Psi_k^{\omega_i} = \lambda_k^{\omega_i} S^{\omega_i}\Psi_k^{\omega_i} ,\quad
A^{+,\omega_i}\Psi_k^{\omega_i} = \lambda_k^{\omega_i} S^{+, \omega_i} \Psi_k^{\omega_i}
\label{offeig4}
\end{equation}
for the $k$-th eigenpair.
We present
two choices: one with no oversampling and one with oversampling,
though one can consider various options \cite{eglp13oversampling}.
In (\ref{offeig4}), with no oversampling,
we can choose $A^{\omega_i}= [a^{\omega_i}_{mn}]$ and $S^{\omega_i} = [s^{\omega_i}_{mn}]$, where
$ a^{\omega_i}_{mn} = \int_{\omega_i} {\kappa} \nabla \psi_m^{\omega_i} \cdot \nabla \psi_n^{\omega_i} = (R^{\omega_i}_{\text{snap}})^T {A} R^{\omega_i}_{\text{snap}}$ and
\quad\text{and}\quad
$s^{\omega_i}_{mn} = \int_{\omega_i} \widetilde{{\kappa}} \,\psi_m^{\omega_i}\, \psi_n^{\omega_i} = (R^{\omega_i}_{\text{snap}})^T {S} R_{\text{snap}}^{\omega_i}$.
With oversampling, we can choose
$A^{+,\omega_i}= [a^{+,\omega_i}_{mn}]$ and $S^{+,\omega_i} = [s^{+,\omega_i}_{mn}]$,
where
$ a^{+,\omega_i}_{mn} = \int_{\omega_i^+} {\kappa} \nabla \psi_m^{+,\omega_i} \cdot \nabla \psi_n^{+,\omega_i} = (R^{+,\omega_i}_{\text{snap}})^T {A} R^{+,\omega_i}_{\text{snap}}$ and
$s^{+,\omega_i}_{mn} = \int_{\omega_i^+} \widetilde{{\kappa}} \,\psi_m^{+,\omega_i}\, \psi_n^{+,\omega_i} = (R^{+,\omega_i}_{\text{snap}})^T {S} R_{\text{snap}}^{+,\omega_i}$.
%$A^{+,\text{off}} = [a_{mn}^{+,\text{off}}] = \int_{\omega_i^+} {\kappa}(x) \nabla \psi_m^{+,\text{snap}} \cdot \nabla \psi_n^{+,\text{snap}} = \big(R_{\text{snap}}^+\big)^T \overline{A}^+ R_{\text{snap}}^+$,
%$S^{+,\text{off}} = [s_{mn}^{+,\text{off}}] = \int_{\omega_i^+} \widetilde{{\kappa}}(x) \psi_m^{+,\text{snap}}  \psi_n^{+,\text{snap}} = \big(R_{\text{snap}}^+\big)^T {S}^+ R_{\text{snap}}^+$.
To generate the offline space we then choose the smallest $l_i$ eigenvalues  and form the corresponding eigenvectors in the respective space of snapshots by setting
$\phi_k^{\omega_i} = \sum_{j=1}^{L_i} \Psi_{kj}^{\omega_i} \psi_j^{\omega_i}$ or $\phi_k^{+,\omega_i} = \sum_{j=1}^{L_i^+} \Psi_{kj}^{+,\omega_i} \psi_j^{+,\omega_i}$
(for $k=1,\ldots, L_i$ or $k=1,\cdots, L_i^+$), where $\Psi_{kj}^{\omega_i}$ and
$\Psi_{kj}^{+,\omega_i}$ are the coordinates of the vector $\Psi_{k}^{\omega_i}$ and $\Psi_{k}^{+,\omega_i}$ respectively.
Collecting all offline basis functions and using a single index notation,
we then create the offline matrices
 $
R_{\text{off}}^+ = \left[ \psi_{1}^{+,\text{off}}, \ldots, \psi_{M_{\text{off}}}^{+,\text{off}} \right]$
and
$R_{\text{off}} = \left[ \psi_{1}^{\text{off}}, \ldots, \psi_{M_{\text{off}}}^{\text{off}} \right]
$
to be used in the online space construction, where $M_{\text{off}}$ is the total number of offline basis functions.
The discrete system corresponding
to (\ref{eq:GMsFEM_global}) is
\[
R_{\text{off}}^T A R_{\text{off}} u_H^{\text{discrete}}= R_{\text{off}}^T f
\]
where $u_H^{\text{discrete}}$ is the discrete version of $u_H$.

\subsection{A numerical example}

\begin{figure}[h]
\centering
\includegraphics[scale=0.4]{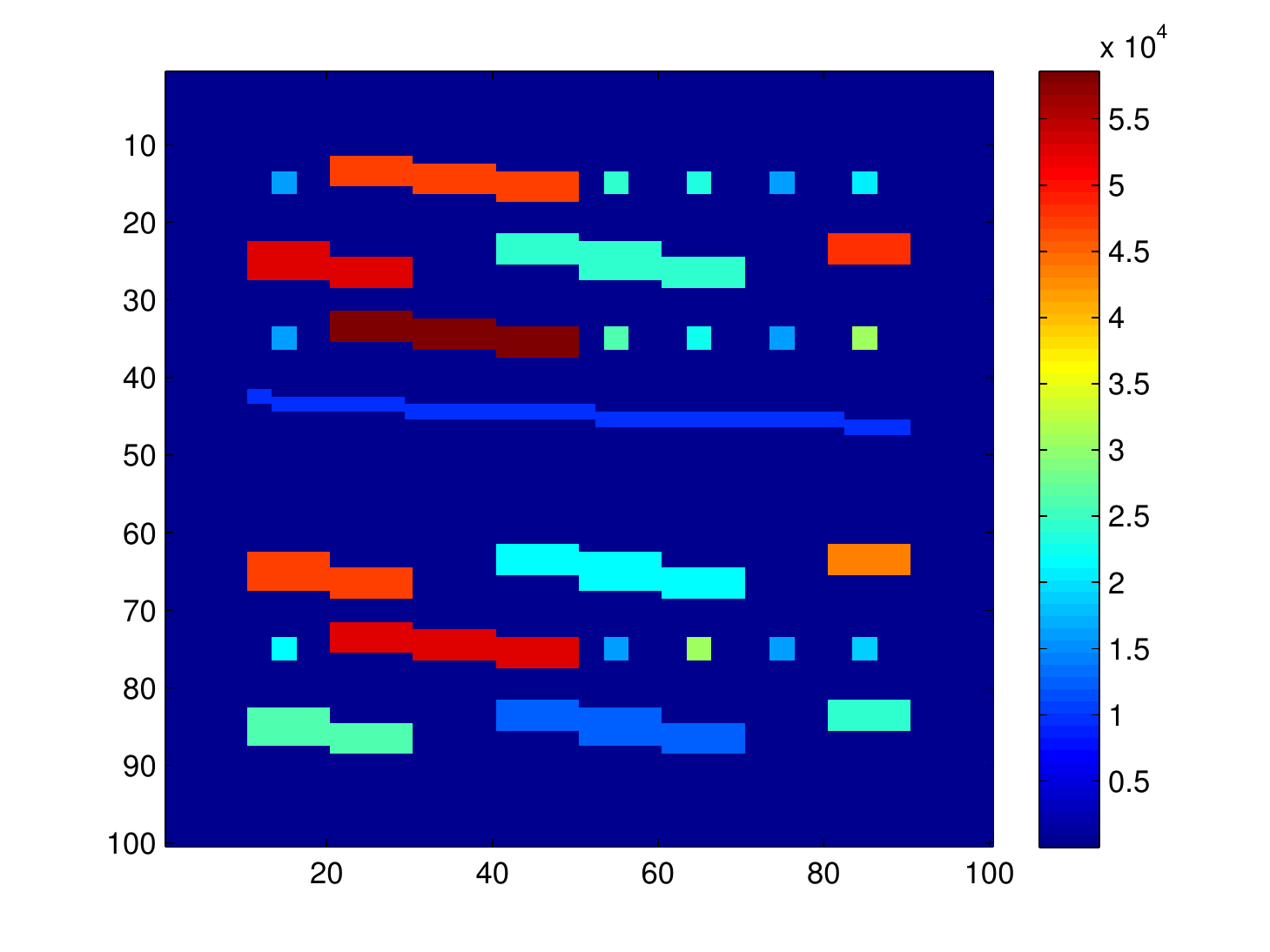}
\caption{The permeability field $\kappa(x)$.}
\label{fig:perm_random}
\end{figure}

We present a numerical result that demonstrate the convergence
of the GMsFEM. More detailed numerical studies can be found
in the literature. We consider
the permeability field, $\kappa(x)$ (see (\ref{main:red}))
that is shown in Figure \ref{fig:perm_random},
the source term $f=0$, and the boundary condition to be $x_1$. The fine grid
is $100\times 100$, coarse grid is $10\times 10$. We consider
the snapshot space spanned by harmonic functions in the oversampled
domain (with randomized boundary conditions)
and vary the number of basis functions per node.
The numerical results are shown in Table \ref{tab:GMsFEM_results}
in which $e_a$ and $e_2$ denote errors in energy and $L^2$ norms repsectively.
As
we observe from these numerical results  the GMsFEM converges as we
increase the number of basis functions. We also show the convergence
when using polynomial basis functions. It is clear that the FEM
method with polynomial basis functions does not converge
and one can only
observe this convergence only after a very large number of
 basis functions are chosen so that we cross the fine-scale threshold.
We note that the convergence with one basis function per node does not
perform well for the GMsFEM because of the high contrast.
The first, second, and third
smallest eigenvalues, $\Lambda_*$, among all coarse blocks are
 $0.0024$, $24.0071$, $35.6941$, respectively.
We emphasize that if the eigenvalue problem
(and the snapshot space) is not chosen appropriately,
one can get more very small eigenvalues (\cite{ge09_2}).
As we see that the first
smallest eigenvalue is very small and, as a result, the error is large,
when using one basis function.
The eigenvalue distribution is important for online basis construction,
as will be discussed later.
%55.1494, 79.1868

\begin{table}[!htb]
\centering
%\begin{tabular}{|c|c|c|}
%\hline
%\#basis (DOF) & $e_{a}$ & $e_{2}$\tabularnewline
%\hline
%1 (81) & 47.02\% & 28.40\%\tabularnewline
%\hline
%2 (162) & 24.66\% & 6.04\%\tabularnewline
%\hline
%3 (243) & 23.10\% & 5.32\%\tabularnewline
%\hline
%4 (324) & 21.51\% & 4.64\%\tabularnewline
%\hline
%5 (405) & 17.41\% & 3.07\%\tabularnewline
%\hline
%\end{tabular}
\begin{tabular}{|c|c|c|}
\hline
\#basis (DOF) & $e_{a}$ & $e_{2}$\tabularnewline
\hline
1 (81) & 69.05\% & 12.19\%\tabularnewline
\hline
2 (162) & 22.55\% & 1.19\%\tabularnewline
\hline
3 (243) & 19.86\% & 0.99\%\tabularnewline
\hline
4 (324) & 16.31\% & 0.70\%\tabularnewline
\hline
5 (405) & 14.20\% & 0.65\%\tabularnewline
\hline
\end{tabular}
\begin{tabular}{|c|c|c|}
\hline
DOF & $e_{a}$ & $e_{2}$\tabularnewline
\hline
81 & 103.2\% & 23\%\tabularnewline
\hline
361 & 100\% & 23\%\tabularnewline
\hline
841 & 80.33\% & 15\%\tabularnewline
\hline
\end{tabular}
\caption{Left table: The convergence for the GMsFEM using multiscale basis functions. Right table: The convergence for FEM using polynomials basis functions.}
\label{tab:GMsFEM_results}
\end{table}

% \begin{table}[h]
% \centering
% \begin{tabular}{|c|c|c|}
% \hline
% \#basis & $e_{a}$ & $e_{2}$\tabularnewline
% \hline
% 1 (81) & 47.02\% & 28.40\%\tabularnewline
% \hline
% 2 (162) & 24.66\% & 6.04\%\tabularnewline
% \hline
% 3 (243) & 23.10\% & 5.32\%\tabularnewline
% \hline
% 4 (324) & 21.51\% & 4.64\%\tabularnewline
% \hline
% 5 (405) & 17.41\% & 3.07\%\tabularnewline
% \hline
% \end{tabular}
% \caption{The convergence for different number of basis functions for the GMsFEM.}
% \label{tab:GMsFEM_results}
% \end{table}

\section{Adaptivity in GMsFEM}
\label{sec:adaptivity}

The success of the GMsFEM depends on adaptive implementation and appropriate
error indicators.
In this section, we discuss an a-posteriori error indicator for the GMsFEM framework.
We will demonstrate the main idea using a
continuous Galerkin formulation; however,
this concept can be generalized to other discretizations as we will discuss
later on.
This error indicator is further used
to develop an adaptive enrichment algorithm for
the linear elliptic equation with multiscale high-contrast coefficients.
Rigorous a-posteriori
error indicators are needed to perform an adaptive enrichment.
We would like to point out that
there are many related activities in designing a-posteriori
error estimates
%\cite{Dorfler96,ohl12, abdul_yun, dinh13, nguyen13, tonn11}
 for global reduced models.
The main difference is that the error estimators presented in this
section
 are
based on a special local eigenvalue problem and use the eigenstructure
of the offline space.

We can consider various kinds of error indicators
that are based on
the $L^2$-norm of the local residual
and the other is based on the weighted $H^{-1}$-norm (we will also call it
$H^{-1}$-norm based) of the local residual
where the weight is related to the coefficient of the elliptic equation.
The latter will be studied in the paper and we refer to
\cite{chung2014adaptive}
for more details.

Let $u_{H} \in V_{H,\text{off}}$ be the solution obtained in
(\ref{eq:GMsFEM_global}).
Consider a given coarse neighborhood $\omega_i$.
We define a space $V_i = H^1_0(\omega_i) \cap V$
which is equipped with the norm $\|v\|_{V_i}^2 =\int_{\omega_i} \kappa(x) | \nabla v|^2 $.
We also define the following linear functional on $V_i$ by
\begin{equation}
R_i(v) =  \int_{\omega_i} fv - \int_{\omega_i} \kappa\nabla u_{H}\cdot \nabla v.
\end{equation}
This is called the $H^{-1}$-residual on $\omega_i$,
The functional norm of $R_i$, denoted by $\|R_i\|_{V_i^*}$,
gives a measure of the size of the residual.
The first important result
% in \cite{Adaptive-GMsFEM}
states that these residuals give a computable
indicator of the error $u-u_{H}$
in the energy norm. We have
\begin{equation}
\label{eq:residualbound}
\| u-u_{H}\|_V^2 \leq C_{\text{err}}  \sum_{i=1}^{N_c} \|R_i\|^2_{V_i^*}  (\lambda^{\omega_i}_{l_i+1})^{-1},
\end{equation}
where $C_{\text{err}}$ is a uniform constant,
and $\lambda^{\omega_i}_{l_i+1}$ denotes the $({l_i+1})$-th eigenvalue for the problem (\ref{offeig4})
in the coarse neighborhood $\omega_i$, and corresponds to the first eigenvector that is not included in the construction of $V_{H,\text{off}}^{\omega_i}$.

%The adaptive enrichment algorithm \cite{Adaptive-GMsFEM}
%is stated as follows.
The above error indicator allows one to construct an adaptive enrichment algorithm.
It is an iterative process, and basis functions are added in each iteration based on the current solution.
We use the index $m\geq 1$ to represent the enrichment level.
At the enrichment level $m$, we use $V_{H,\text{off}}^m$
to denote the corresponding GMsFEM space
and $u_{H}^m$ the corresponding solution obtained in (\ref{eq:GMsFEM_global}) using the space
$V_{H,\text{off}}^m$.
Furthermore, we use $l_i^m$ to denote the number of basis functions used in the coarse neighborhood $\omega_i$.
We will present the strategy for getting the space
$V_{H,\text{off}}^{m+1}$ from $V_{H,\text{off}}^m$.
Let $0< \theta < 1$ be a given number independent of $m$.
First of all, we compute the local
residuals for every coarse neighborhood $\omega_i$:
\begin{equation*}
\eta_i^2 = \|R_i\|^2_{V_i^*}  (\lambda^{\omega_i}_{l^m_i+1})^{-1}, \quad\quad i=1,2,\cdots, N_c,
\end{equation*}
where $R_i(v)$ is defined using $u_{H}^m$, namely,
\begin{equation*}
R_i(v) =  \int_{\omega_i} fv - \int_{\omega_i} \kappa\nabla u^m_{H}\cdot \nabla v, \quad \forall v\in V_i.
\end{equation*}
Next, we will add basis functions for the coarse neighborhoods with large residuals.
To do so, we re-enumerate the coarse neighborhoods so that the above local residuals $\eta_i^2$ are arranged in decreasing order
$\eta^2_1 \geq \eta^2_2 \geq \cdots \geq \eta^2_N$.
We then select the smallest integer $k$ such that
\begin{equation}
\label{eq:criteria_online}
\theta \sum_{i=1}^{N_c} \eta_i^2 \leq \sum_{i=1}^k \eta_i^2.
\end{equation}
For those coarse neighborhoods $\omega_1, \cdots, \omega_k$ (in the new enumeration)
chosen in the above procedure, we will add basis functions
by using the next eigenfunctions $\Psi_{l_i+1}^{\omega_i}, \Psi_{l_i+2}^{\omega_i}, \cdots$.
We remark that the number of new additional basis depends on the eigenvalue decay.
The resulting space is called $V_{\text{off}}^{m+1}$.
We remark that the choice of $k$ defined in (\ref{eq:criteria_online})
is called the Dorlfer's bulk marking strategy \cite{Dorfler96}.

As we mentioned
the local basis functions do not contain any global information, and thus they
cannot be used to efficiently
capture these global behaviors.
We will therefore present in the next section
a fundamental ingredient of GMsFEM - the development of
online basis functions that is
necessary to obtain a coarse representation of the fine-scale solution
and
gives a rapid convergence of the corresponding adaptive enrichment algorithm.

\begin{remark}[Implementation]
The algorithm above can be described as follows.
We start with an initial space with a small number of basis functions for each coarse grid block. Then we solve the problem and compute the error estimator. We locate the coarse grid blocks with large errors and add more basis functions for these coarse grid blocks. This procedure is repeated until the error goes below a certain tolerance. We remark that the adaptive strategy belongs to the online process, because it is the actual simulation.
On the other hand, the generation of basis functions belongs to the offline process. About stopping criteria for this algorithm, one can stop the algorithm when the total number of basis functions reach a certain level.
On the other hand, one can stop the algorithm when the value of the error indicator goes below a certain tolerance.
\end{remark}

\begin{remark}[Goal Oriented Adaptivity]
For some practical problems, one is interested in approximating some
function of the solution, known as the quantity of interest, rather than
the solution itself.  Examples include an average or weighted
average of the solution over a particular
subdomain, or some localized solution response. In these cases,
goal-oriented adaptive methods yield a more efficient approximation
than standard adaptivity, as the enrichment of degrees of freedom is
focused on the local improvement of the quantity of interest rather than
across the entire solution. In \cite{chung2016goal},
 we study goal-oriented adaptivity for multiscale methods, and
in particular the design of error indicators to drive the adaptive enrichment
based on the goal function.
In this methodology, one seeks to determine the number of multiscale basis
functions adaptively for each coarse region
to efficiently reduce the error in the goal functional.
Two estimators are studied. One is a residual based strategy
and the other uses dual weighted residual method for multiscale problems.
The method is demonstrated on high-contrast problems with
heterogeneous multiscale coefficients, and is seen to outperform
the standard residual based strategy with respect to efficient reduction
of error in the goal function.

\end{remark}

%************************************
\subsection{Numerical Results}
\label{sec:numresults_adaptive}
%************************************

\begin{figure}[h]
\centering
\includegraphics[scale=0.4]{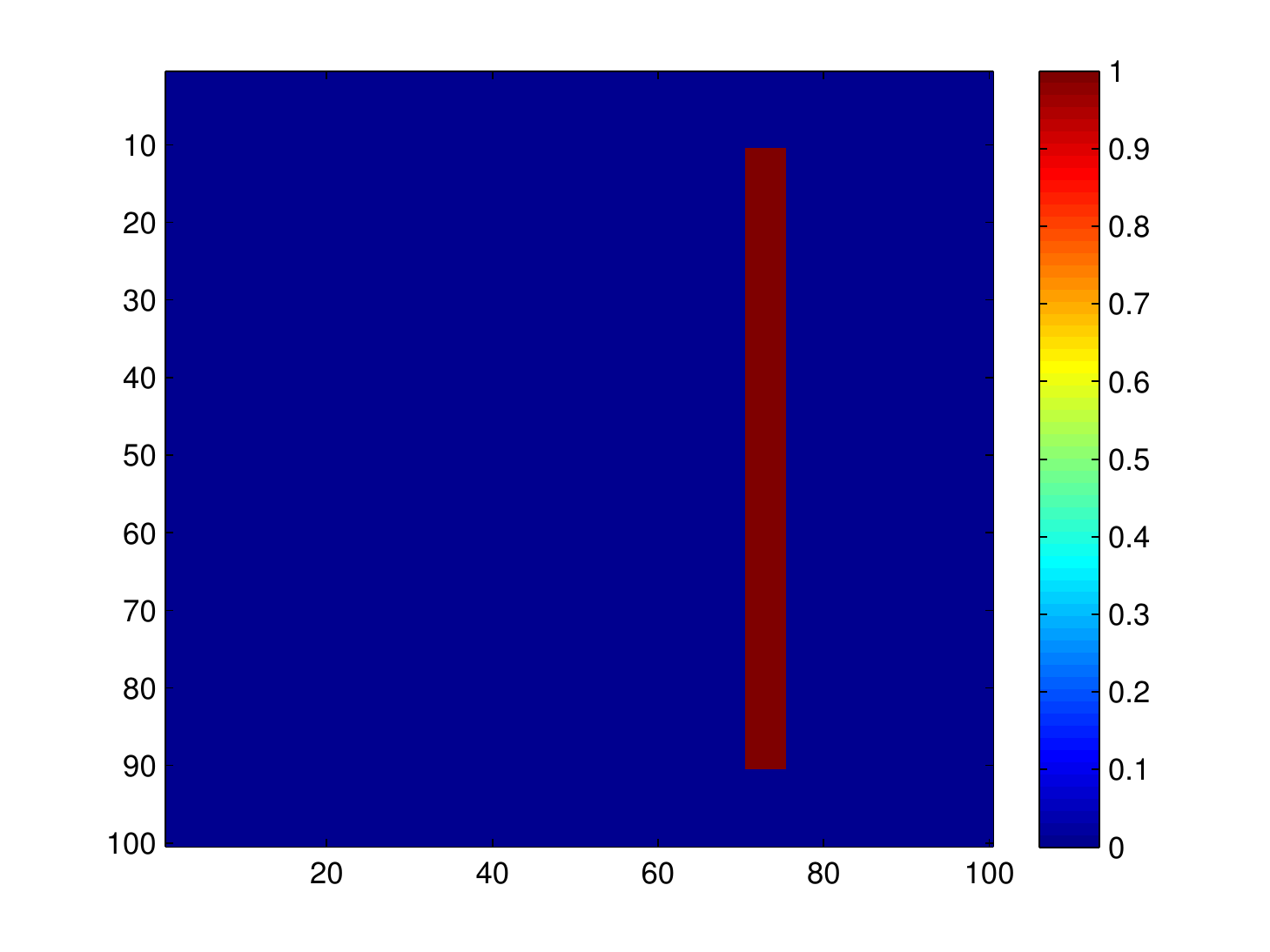}
\includegraphics[scale=0.4]{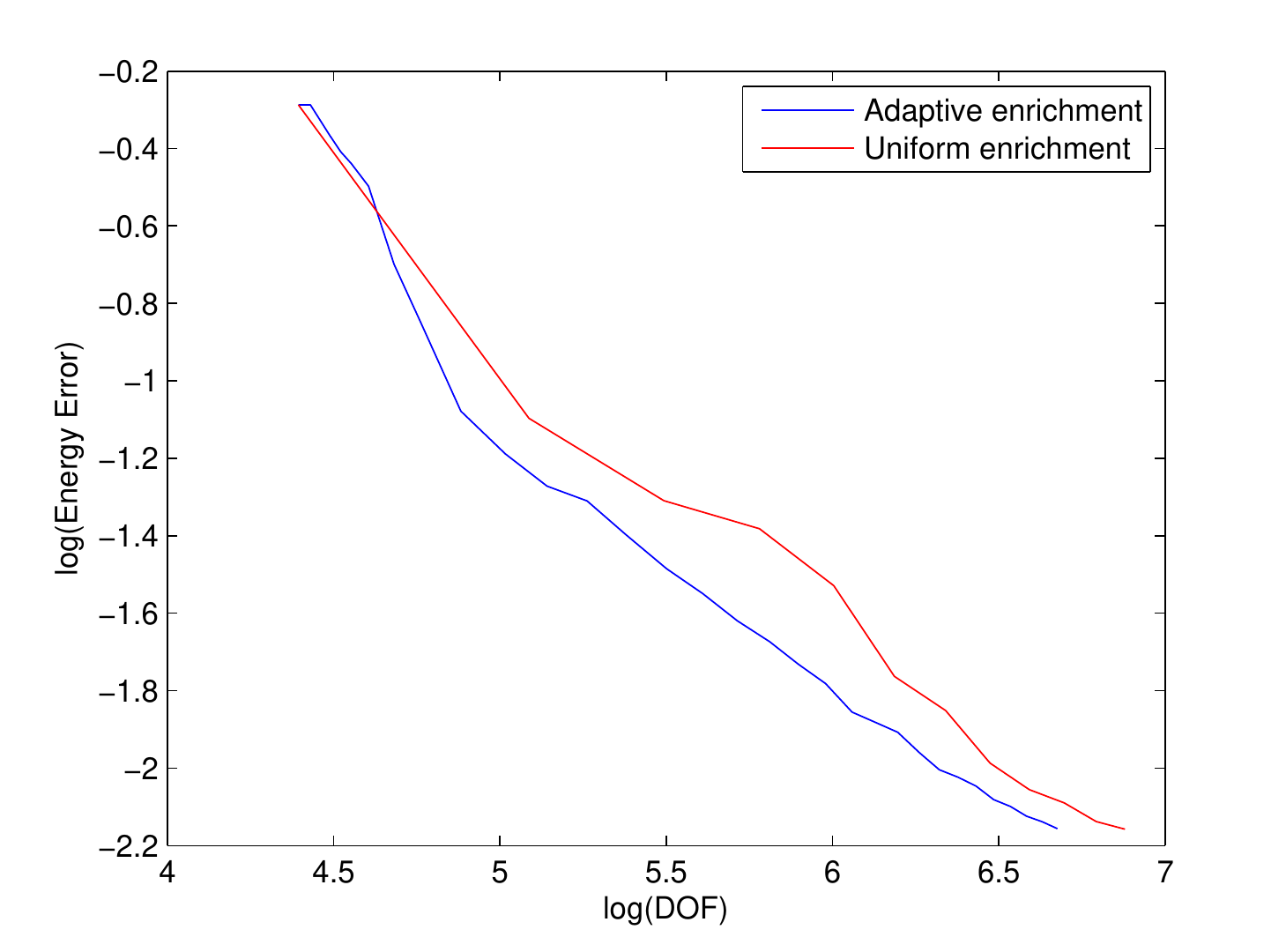}
\protect\caption{Left: Source function $f$. Right: The convergence of adaptive vs. uniform enrichment.}
\label{fig:source_term}
\end{figure}

We show numerical results for the $H^{-1}$ adaptivity. We consider
the permeability field shown in Figure \ref{fig:perm_random}
for (\ref{main:red})
and the forcing term
shown in Figure \ref{fig:source_term}. In Table \ref{tab:GMsFEM_adaptive},
we present the numerical results and in Figure \ref{fig:source_term}
(right plot), we compare the convergence of the adaptive GMsFEM and the GMsFEM,
which uses a uniform number of basis functions. As we observe that
the adaptive GMsFEM converges faster.

%
%\begin{table}[h]
%%\label{tab:offline_adaptive_results}
%\centering
%\begin{tabular}{|c|c|c|}
%\hline
%\#DOF & $e_{a}$ & $e_{2}$\tabularnewline
%\hline
%81 & 75.05\% & 42.48\%\tabularnewline
%\hline
%116 & 42.66\% & 15.21\%\tabularnewline
%\hline
%268 & 20.04\% & 3.84\%\tabularnewline
%\hline
%487 & 11.56\% & 1.52\%\tabularnewline
%\hline
%727 & 8.00\% & 0.83\%\tabularnewline
%\hline
%\end{tabular} %
%\begin{tabular}{|c|c|c|}
%\hline
%\#DOF & $e_{a}$ & $e_{2}$\tabularnewline
%\hline
%81 & 75.05\% & 42.48\%\tabularnewline
%\hline
%162 & 32.71\% & 6.30\%\tabularnewline
%\hline
%243 & 24.86\% & 5.04\%\tabularnewline
%\hline
%486 & 15.90\% & 2.00\%\tabularnewline
%\hline
%729 & 10.94\% & 1.10\%\tabularnewline
%\hline
%\end{tabular}
%\protect\caption{Errors for the adaptive GMsFEM.
% Left: Adaptive enrichment. Right: Uniform enrichment. }
%\label{tab:GMsFEM_adaptive}
%\end{table}

\begin{table}[h]
%\label{tab:offline_adaptive_results}
\centering
\begin{tabular}{|c|c|c|}
\hline
\#DOF & $e_{a}$ & $e_{2}$\tabularnewline
\hline
81 & 75.04\% & 42.48\%\tabularnewline
\hline
151 & 30.47\% & 7.84\%\tabularnewline
\hline
245 & 22.65\% & 4.70\%\tabularnewline
\hline
334 & 18.76\% & 3.59\%\tabularnewline
\hline
395 & 16.84\% & 3.08\%\tabularnewline
\hline
\end{tabular} %
\begin{tabular}{|c|c|c|}
\hline
\#DOF & $e_{a}$ & $e_{2}$\tabularnewline
\hline
81 & 75.04\% & 42.48\%\tabularnewline
\hline
162 & 33.39\% & 6.74\%\tabularnewline
\hline
243 & 27.00\% & 5.52\%\tabularnewline
\hline
324 & 25.11\% & 4.72\%\tabularnewline
\hline
405 & 21.68\% & 3.50\%\tabularnewline
\hline
\end{tabular}
\protect\caption{Errors for the adaptive GMsFEM.
 Left: Adaptive enrichment. Right: Uniform enrichment. }
\label{tab:GMsFEM_adaptive}
\end{table}

\section{Residual-based online procedure}
 \label{sec:online}

Previously, we discussed
adaptive enrichment procedures and derived an
a-posteriori error indicator,
which gives an estimate of the local error
over coarse grid regions. We developed the error indicators based
 on
the $L^2$-norm of the local residual
and on the weighted $H^{-1}$-norm of the local residual,
where the weight is related to the coefficient of the elliptic equation.
%We have shown that the use of weighted $H^{-1}$-norm residual
%gives a more robust error indicator
%which works well for cases with high contrast media.
Adaptivity is important for local multiscale methods as it  identifies
regions with large errors. However, after
adding some initial basis functions,
one needs to take into account some global information as the distant effects
can be important. In this section, we discuss the development of online
basis functions that substantially accelerate the convergence of GMsFEM.
The online basis functions are constructed based on a residual and motivated
by the analysis.
We show
 that one needs to have a sufficient
number of initial (offline) basis functions   to guarantee
an error decay independent of the contrast.
{ In particular, the error decay in one adaptive iteration
 is proportional to $1-C\Lambda_*$ (see (\ref{eq:error_decay_online}) for
more precise estimate), where  $\Lambda_*$
is the smallest eigenvalue among all coarse blocks that
the corresponding eigenvector is not included in the offline space.
Thus, it is important to include all eigenvectors corresponding to very small
eigenvalues in the offline space.
}
In general, we would like to apply one iteration of the online procedure
and reach a desired error threshold.
Numerical results are presented to demonstrate
that one needs to have a sufficient number of initial basis
functions in the offline space before constructing
online multiscale basis functions.

\subsection{Residual-based online adaptive GMsFEM}

We use the index $m\geq 1$ to represent the enrichment level.
At the enrichment level $m$, we use $V_{H,\text{off}}^m$
to denote the corresponding GMsFEM space
and $u_{H}^m$ the corresponding solution.
The sequence of functions $\{ u_{H}^m \}_{m \geq 1}$ will converge
to the fine-scale solution.
We emphasize that the space $V_{H,\text{off}}^m$
can contain both offline and online basis functions.
We will construct a strategy for getting the space
$V_{H,\text{off}}^{m+1}$ from $V_{H,\text{off}}^m$.
The online basis functions are computed based on some local residuals
for the current multiscale solution, that is, the function $u_{H,\text{off}}^m$.

Suppose that we need to add a basis function $\phi^{\text{on},\omega_i} \in V_i$
on the $i$-th coarse neighbourhood $\omega_i$.
%where $V_h^{\omega_i}$ is the fine-scale space restricted in $\omega_i$ with zero boundary conditions.
Let $V_{H,\text{off}}^{m+1} = V_{H,\text{off}}^m + \text{span} \{ \phi^{\text{on},\omega_i} \}$ be the new approximation space,
and $u_{H}^{m+1} \in V_{H,\text{off}}^{m+1} $ be the corresponding GMsFEM solution.
It is easy to see  that $u_{H}^{m+1}$ satisfies
\begin{equation*}
\| u - u_{H}^{m+1}\|_{V}^2 = \inf_{v \in V_{H,\text{off}}^{m+1}}  \| u - v\|_V^2.
\end{equation*}
Taking $v = u_{H}^{m+1} + \alpha \phi^{\text{on},\omega_i}$,
where $\alpha$ is a scalar to be determined, we have
\begin{equation}
\label{eq:error_decay_online1}
\| u - u_{H}^{m+1} \|_{V}^2 \leq \| u - u_{H}^{m} - \alpha \phi^{\text{on},\omega_i} \|_V^2
= \|u- u_{H}^{m}\|^2_V - 2 \alpha a(u- u_{H}^{m}, \phi^{\text{on},\omega_i}) + \alpha^2 a(\phi^{\text{on},\omega_i},\phi^{\text{on},\omega_i}).
\end{equation}
The last two terms in the above inequality measure the amount of reduction in error
when the new basis function $\phi^{\text{on},\omega_i}$
is added to the space $V_{H,\text{off}}^{m} $.
To determine $\phi^{\text{on},\omega_i}$,
we first assume that the basis function $\phi^{\text{on},\omega_i}$
is normalized so that $a(\phi^{\text{on},\omega_i},\phi^{\text{on},\omega_i})=1$.
In order to maximize the reduction in error, one can show that
$\alpha = a(u-u_{H}^{m},\phi^{\text{on},\omega_i})$.
Using this choice of $\alpha$, we have
\begin{equation*}
\| u - u^{m+1}_{H}\|_{V}^2 \leq \|u-u^m_{H}\|^2_V - |a(u-u^m_{H},\phi^{\text{on},\omega_i})|^2.
\end{equation*}
%Since $\phi \in V_i \subset V$, by using (\ref{eq:fineprob}), we have
%\begin{equation*}
%\| u - u^{m+1}_{\text{ms}}\|_{V}^2 \leq \|u-u^m_{\text{ms}}\|^2_V - |(f,\phi)-a(u^m_%{\text{ms}},\phi)|^2.
%\end{equation*}
%We will then find $\phi\in V_i$ to maximize the local residual $|(f,\phi)-a(u^m_{\te%xt{ms}},\phi)|^2$.
%Clearly, the maximum of the quantity $|(f,\phi)-a(u^m_{\text{ms}},\phi)|$
%equals to the functional norm of the residual $R_i$.
%Moreover,
%the required $\phi \in V_i$ is the solution of
We can show that
\begin{equation}
\label{eq:online}
a(\phi^{\text{on},\omega_i},v) = (f,v)-a(u^m_{H},v), \quad \forall v\in V_i
\end{equation}
and $\|\phi^{\text{on},\omega_i}\|_{V_i} = \|R_i\|_{V_i^*}$.
Hence, the new online basis function
$\phi^{\text{on},\omega_i}\in V_i$ can be obtained by solving (\ref{eq:online}).
In addition, the residual norm $\|R_i\|_{V_i^*}$
provides a measure on the amount of reduction in energy error.

In \cite{chung2015residual, chung2015online},
we have studied the convergence of
the above online adaptive procedure.
To simplify notations, we write $r_i = \| R_i \|_{V_i^*}$.
%From the above constructions, we have
%\begin{equation}
%\| u - u^{m+1}_{\text{ms}}\|_{V}^2 \leq \|u-u^m_{\text{ms}}\|^2_V - r_i^2.
%\label{eq:err1_online}
%\end{equation}
%We assume that each of the spaces $V_{\text{ms}}^m$, $m\geq 1$,
%contains $n_j$ offline basis functions for the coarse neighborhood $\omega_j$.
%Then, similar to (\ref{eq:residualbound}), we have
%\begin{equation}
%\| u - u^m_{\text{ms}}\|_V^2 \leq C_{\text{err}}  \sum_{j=1}^N r_j^2 (\lambda_{n_j+1%}^{\omega_j})^{-1}.
%\label{eq:err2_online}
%\end{equation}
%Combining (\ref{eq:err1_online}) and (\ref{eq:err2_online}), we obtain
We have shown that
\begin{equation}
\label{eq:error_decay_online}
\| u - u^{m+1}_{H}\|_{V}^2
\leq \Big(1 - \frac{\lambda_{n_i+1}^{\omega_i}}{C_{\text{err}}} \frac{ r_i^2 (\lambda_{n_i+1}^{\omega_i})^{-1} }{  \sum_{j=1}^{N_c} r_j^2 (\lambda_{n_j+1}^{\omega_j})^{-1} } \Big) \| u - u^m_{H}\|_V^2,
\end{equation}
where $n_i$ is the number of offline basis in $\omega_i$.
The above inequality gives the convergence of the online adaptive GMsFEM
with a precise convergence rate for the case when one online
basis function is added per iteration.
{ The estimate (\ref{eq:error_decay_online}) shows that the eigenvectors
corresponding to very small eigenvalues need to be included in the offline
space in order to achieve a significant error reduction in one iteration.}
To enhance the convergence and efficiency of the
online adaptive GMsFEM, we consider enrichment on non-overlapping coarse neighborhoods.
Let $I \subset \{ 1,2,\cdots, N\}$ be the index set of some non-overlapping coarse neighborhoods.
For each $i\in I$, we can obtain a basis function $\phi^{\text{on},\omega_i}\in V_i$ using (\ref{eq:online}).
We define $V_{H,\text{off}}^{m+1} = V_{H,\text{off}}^m + \text{span} \{ \phi^{\text{on},\omega_i} \, , \, i\in I\}$.

We
remark that the error will decrease
independent of physical  parameters such as the contrast and scales
if the offline space is appropriately chosen.
We will demonstrate the effectiveness of this method by a numerical example.
We remark that one can also derive apriori error estimate for
$\| u - u^{m+1}_{H,\text{off}}\|_{V}^2$.

\subsection{Numerical result}
\label{sec:numresults_online}

In Table \ref{tab:GMsFEM_online}, we present numerical
results for online enrichment for the GMsFEM. We use the permeability
field shown in Figure \ref{fig:perm_random} and the forcing
term shown in Figure \ref{fig:source_term}
(left figure). As we observe
the online enrichment does not improve
the error if we have only one offline basis function.
We remind that
the first, second, and third
smallest eigenvalues, $\Lambda_*$, among all coarse blocks are
 $0.0024$, $24.0071$, $35.6941$, respectively. Because the first eigenvalue
is small, the error decrease in one online iteration is small.
 In particular, for each online iteration, the error
decreases slightly. It is important that one adaptive online iteration
can decrease the error substantially.
As we increase the number of offline basis functions, the convergence
is very fast and one online iteration is sufficient to reduce
the error significantly. In one iteration, the error drops below $1$\%
for the energy error. We note that this procedure can be implemented
adaptively and we add online basis functions only in some regions.

%\begin{table}[h]
%\centering
%\begin{tabular}{|c|c|c|c|c|c|c|c|}
%\hline
%DOF\textbackslash{}\#initial basis & 1 & 2 & 3 &  & 1 & 2 & 3\tabularnewline
%\hline
%81 & 75.06\% & - & - &  & 42.49\% & - & -\tabularnewline
%\hline
%162 & 32.77\% & 32.76\% & - &  & 16.29\% & 6.32\% & -\tabularnewline
%\hline
%243 & 21.59\% & 0.89\% & 24.94\% &  & 7.33\% & 0.060\% & 5.07\tabularnewline
%\hline
%324 & 3.45\% & 0.013\% & 1.06\% &  & 1.21\% & 0.0011\% & 0.064\tabularnewline
%\hline
%405 & 2.45\% & 1.62e-04\% & 0.13\% &  & 0.68\% & 1.31e-05\% & 7.08e-04\tabularnewline
%\hline
%486 & 1.18\% & 1.32e-06\% & 1.72e-04\% &  & 0.31\% & 9.00e-07\% & 1.35e-05\%\tabularnewline
%\hline
%\end{tabular}
%\protect\caption{The errors for the GMsFEM with the online procedure (left
%table: $H^1$ errors; right table: $L^2$ errors).}
%\label{tab:GMsFEM_online}
%\end{table}

\begin{table}[h]
\centering
\begin{tabular}{|c|c|c|c|c|c|c|c|}
\hline
DOF\textbackslash{}\#initial basis & 1 & 2 & 3 &  & 1 & 2 & 3\tabularnewline
\hline
81 & 75.06\% & - & - &  & 42.49\% & - & -\tabularnewline
\hline
162 & 32.77\% & 33.50\% & - &  & 16.29\% & 6.76\% & -\tabularnewline
\hline
243 & 21.59\% & 1.13\% & 27.14\% &  & 7.33\% & 0.079\% & 5.54\tabularnewline
\hline
324 & 3.46\% & 0.019\% & 1.35\% &  & 1.21\% & 0.0018\% & 0.081\tabularnewline
\hline
405 & 2.45\% & 2.61e-04\% & 0.018\% &  & 0.68\% & 1.81e-05\% & 0.0011\tabularnewline
\hline
486 & 1.18\% & 2.88e-06\% & 2.49e-04\% &  & 0.31\% & 1.74e-07\% & 1.60e-05\%\tabularnewline
\hline
\end{tabular}
\protect\caption{The errors for the GMsFEM with the online procedure (left
table: $H^1$ errors; right table: $L^2$ errors).}
\label{tab:GMsFEM_online}
\end{table}

\section{Selected global discretizations and energy minimizing oversampling}
\label{sec:selected_discretization}

In previous sections, we used continuous Galerkin coupling for multiscale
basis functions. In many applications, one needs to use various
discretizations. For example, for flows in porous media, mass conservation
is very important and, thus, it is advantageous to use mixed methods.
  In seismic wave applications involving
the time-explicit discretization of wave equations,
one needs block-diagonal mass matrices, which can be obtained using
discontinuous Galerkin approaches. In this section, we present
two global couplings, mixed finite element and
Interior Penalty Discontinuous
Galerkin (IPDG) couplings. We refer to \cite{efendiev2015spectral} for
 Hybridized Discontinuous Galerkin (HDG) coupling and
to \cite{chung2015_EMO} for comparisons between using
IPDG and HDG couplings with multiscale basis functions.

For each discretization, we define snapshot spaces and local spectral
decompositions in the snapshot spaces. As we emphasized earlier
the construction of
the snapshot spaces and local multiscale spaces depends on the global
discretization. For example, for mixed GMsFEM, the multiscale
spaces are constructed for the velocity field using two neighboring
coarse elements.
We will present some ingredients of GMsFEM introduced earlier in the
construction of multiscale basis functions, e.g.,
oversampling techniques and online multiscale basis
functions.
We also mention a new ingredient for multiscale basis
construction - energy minimizing snapshots.

In a mixed formulation, the problem can be formulated as
\begin{equation}\label{mixed}
    \kappa^{-1}v + \nabla u = 0 \ \mbox{in } \Omega,\ \ \
    \mbox{div}(v) = f \ \mbox{in } \Omega,
\end{equation}
with Neumann boundary condition $v \cdot n = v_\Omega$ on $\partial \Omega$.
%In the primal formulation, the single-phase problem can be formulated as
%\begin{equation}
%\label{primal}
%-\mbox{div} ( \kappa \nabla u ) = f, \quad\mbox{in } D.
%\end{equation}
Depending on a discretization technique, the coarse grid configuration will
vary. For a mixed formulation, coarse grid blocks that share
common edges (faces) will be used in constructing multiscale basis functions.
In a discontinuous Galerkin formulation, the support of multiscale
basis functions is limited to coarse blocks.

\subsection{Discretizations}

\subsubsection{The mixed GMsFEM}
\label{sec:method}

For the mixed GMsFEM,
we will construct basis functions whose supports are
$\omega_E$, which are the two coarse elements
that share a common edge (face) $E$.
In particular, we let
$\mathcal{E}^{H}$ be the set of all edges (faces) of the coarse grid
and let $N_{e}$ be the total number of edges (faces) of
the coarse grid. We define the coarse grid neighborhood
$\omega_{E}$ of a face $E \in \mathcal{E}^{H}$ as
$$\omega_{E} = \bigcup \{K \in \mathcal{T}^{H} : E \in \partial K \}, \quad i=1,2,\cdots, N_e,$$
which is a union of two coarse grid blocks if $E_i$ is an interior edge (face)
(see Figure \ref{schematic_intro}).
For a coarse edge $E_i$, we write $\omega_{E_i} = \omega_i$ to unify the notations.
%Furthermore,
%we denote by $\omega_i^+$ an oversampling region for $\omega_i$ (see Figure \ref{schematic_intro}).
%One can also
%take the oversampled region $\omega_i^+$ with an arbitrary geometry.
%For the convenience of our analysis presented in the next section,
%we will enlarge $\omega_i$ by one coarse element to form $\omega_i^+$ (see Figure \ref{fig:ill_mixed}).%

Next, we define the notations for the solution spaces for the pressure $u$ and the velocity $v$. Let $Q_{H,\text{off}}$ be the space of functions which are constant on each coarse grid block. We will use this space to approximate $u$. For the multiscale approximation of the velocity $v$, we will follow the general framework of the GMsFEM
and construct a multiscale space $V_{H,\text{off}}$ for the velocity.
%We will first construct a set of basis functions $\beta^{(i)}_{\text{snap}}$ supported in the coarse grid neighborhood $\omega_{i}$. We call
%the span of all these basis function
%$V_{H,\text{snap}} = \bigoplus_{E_{i} \in \mathcal{E}^{H}} V_{H,\text{snap}}^{\omega_i}$ the snapshot space, where $V_{H,\text{snap}}^{\omega_i} =  \text{span} \left( \beta^{(i)}_{\text{snap}} \right)$ is the local snapshot space for the coarse neighborhood $\omega_i$. The snapshot space is an extensive set of functions which can be used to approximate the solution $v$. We reduce the snapshot space to a smaller one before we solve the equation. From each $V^{\omega_i}_{H,\text{snap}}$, we select a subset of basis functions $\beta^{(i)}_{\text{ms}}$ by using an appropriate spectral problem to form a smaller dimensional offline space. We denote the local offline space corresponding to $\omega_i$ by $V_{H,\text{off}}^{\omega_i} = \mbox{span} \left(\beta_{\text{ms}}^{(i)} \right)$ and the global offline space as
%$V_{H,\text{off}} = \bigoplus_{E_{i} \in \mathcal{E}^{H}} V^{\omega_i}_{H,\text{off}}$. The size of $V_{H,\text{off}}$ is generally much smaller than $V_{H,\text{snap}}$, but still contains most important features about the heterogeneous coefficient $\kappa$. We will use the space $V_{H,\text{off}}$ to approximate the velocity $v$.
Using the pressure space $Q_{H,\text{off}}$ and the velocity space $V_{H,\text{off}}$, we solve for $u_{H} \in Q_{H,\text{off}}$ and $v_{H} \in V_{H,\text{off}}$ such that
\begin{equation}\label{MsEquation}
\int_{\Omega} \kappa^{-1} v_{H} \cdot w - \int_{\Omega} \mbox{div}(w) u_{H}= 0, \ \forall w \in \dot{V}_{H,\text{off}}, \quad
\int_{\Omega} \mbox{div}(v_{H}) q = \int_{\Omega} fq, \ \forall q \in Q_{H,\text{off}},
\end{equation}
with boundary condition $v_{H} \cdot n = v_{\Omega,H}$ on $\partial \Omega$, where $\dot{V}_{H,\text{off}} = \{ v \in V_{H,\text{off}} : v \cdot n = 0 \mbox{ on } \partial \Omega\}$, and $v_{\Omega,H}$ is the projection of $v_\Omega$ in the multiscale space.
%in the sense that
%$\int_{E_{i}} (v_{\Omega,H} - v_\Omega) \phi \cdot n = 0$, $\forall \phi \in \beta^{(i)}_{\text{snap}} \mbox{ and } E_{i} \subseteq \partial \Omega$,
%and $v_{\Omega,H}$ is constant on each fine grid face.
We remark that we will define the snapshot and offline spaces
so that the functions in $V_{H,\text{off}}$ are globally $H(\text{div})$-conforming, and as a result the normal components of the basis
are continuous across coarse grid edges.
%The resulting scheme gives a smaller, yet rich enough, dimensional snapshot space
%and a more accurate representation of the solution by multiscale basis.

%\begin{figure}[!ht]
%\centering
%\includegraphics[scale=0.5]{Oversample_Edge_grid}
%\protect\caption{A coarse neighborhood (red) and its corresponding oversampling% region (green) for mixed GMsFEM and GMsHDG. An example of fine grid partition %(blue)
%for a coarse element is also shown.}
%\label{fig:ill_mixed}
%\end{figure}

\subsubsection{DG GMsFEM (GMsDGM)}
\label{sec:dg_gmsfem}
For the GMsDGM, the general methodology for
the construction of multiscale basis functions
is similar to the above mixed GMsFEM.
We will construct multiscale basis functions for the approximation
of $u$ in (\ref{main:red}).
The main difference is that the functions in the snapshot space $V_{H,\text{snap}}$
and the offline space $V_{H,\text{off}}$ are supported
in coarse element $K$, instead of the coarse neighborhood $\omega_i$.
Moreover, in the oversampling approach, the oversampled regions $K^+$
are defined by enlarging coarse elements $K$
(see Figure \ref{schematic_intro}).
The  construction of basis functions will be given in the next section.
When the offline space $V_{H,\text{off}}$
is available, we can find the solution $u_{H}\in V_{H,\text{off}}$
such that (see \cite{chung2014adaptiveDG,chung2015online})
\begin{equation}
a_{\text{DG}}(u_{H},q)=(f,q),\quad\forall q\in V_{H,\text{off}},\label{eq:ipdg}
\end{equation}
where the bilinear form $a_{\text{DG}}$ is defined as
\begin{equation}
\begin{split}
a_{\text{DG}}(u,q)=a_{H}(u,q)-\sum_{E\in\mathcal{E}^{H}}\int_{E}\Big(\average{{\kappa}\nabla{u}\cdot{n}_{E}}\jump{q}+\average{{\kappa}\nabla{q}\cdot{n}_{E}}\jump{u}\Big)+
\sum_{E\in\mathcal{E}^{H}}\frac{\gamma}{h}\int_{E}\overline{\kappa}\jump{u} \jump{q} \label{eq:bilinear-ipdg}
\end{split}
\end{equation}
with
$a_{H}({p},{q})=\sum_{K\in\mathcal{T}^{H}}a_{H}^{K}(p,q)$, $
a_{H}^{K}(p,q)=\int_{K}\kappa\nabla p\cdot\nabla q$,
where $\gamma>0$ is a penalty parameter, ${n}_{E}$ is a fixed unit
normal vector defined on the coarse edge $E \in \mathcal{E}^H$.
Note that, in (\ref{eq:bilinear-ipdg}),
the average and the jump operators are defined in the classical way.
Specifically, consider an interior coarse edge $E\in\mathcal{E}^{H}$
and let $K_{L}$ and $K_{R}$ be the two coarse grid blocks sharing
the edge $E$. For a piecewise smooth function $G$, we define
$\average{G}=\frac{1}{2}(G_{R}+G_{L}),\;\jump{G}=G_{R}-G_{L},\;\text{on}\, E$,
where $G_{R}=G|_{K_{R}}$ and $G_{L}=G|_{K_{L}}$ and we assume that
the normal vector ${n}_{E}$ is pointing from $K_{R}$ to $K_{L}$.
Moreover, on the edge $E$, we define $\overline{\kappa} = (\kappa_{K_R}+\kappa_{K_L})/2$
where $\kappa_{K_{L}}$ is the maximum value of $\kappa$ over $K_{L}$
and $\kappa_{K_R}$ is defined similarly.
For a coarse edge $E$ lying on the boundary $\partial D$, we define
$\average{G}=\jump{G}=G,\; \text{and}\; \overline{\kappa} = \kappa_{K} \;\text{on}\, E$,
where we always assume that ${n}_{E}$ is pointing outside of $D$.
%For vector-valued functions, the above average and jump operators
%are defined component-wise.
We note that the DG coupling (\ref{eq:ipdg})
is the classical interior penalty discontinuous Galerkin (IPDG) method
%\cite{IPDGbook}
with the multiscale basis functions as the approximation space.

%\begin{figure}[!ht]
%\centering
%\includegraphics[scale=0.5]{Oversample_DG_grid_1}
%\protect\caption{A coarse element (red) and its corresponding oversampling region (green) for GMsDGM. An example of fine grid partition (blue)
%for a coarse element is also shown.}
%\label{fig:ill_DG}
%\end{figure}

\subsection{Basis construction}

\subsubsection{Multiscale basis functions in mixed GMsFEM. Non-oversampling}
\label{sec:mixed_emo}

First, we
present the construction of the snapshot space for the approximation
of the velocity $v$.
The space is a large function space containing basis functions
whose normal traces on coarse grid edges
are resolved up to the fine grid level.
Let $E_{i} \in \mathcal{E}^{H}$ be a coarse edge.
We can write the edge
$E_{i} = \bigcup^{J_i}_{j=1} e_{j}$,
where the $e_{j}$'s are the fine grid edges (faces) contained in
$E_i$ and $J_{i}$ is the total number of those fine grid edges.
To construct the snapshot vectors, for each $j=1,2,\cdots, J_i$,
we will solve the following local problem
\begin{equation} \label{localEquation_coupling}
\kappa^{-1} \psi_{i,j} + \nabla \eta_{i,j} = 0 \ \mbox{in } K\subset \omega_{i}, \ \ \
\mbox{div} (\psi_{i,j} ) = \alpha_{i,j} \ \mbox{in } K \subset \omega_{i} ,
\end{equation}
subject to the Neumann boundary conditions $\psi_{i,j} \cdot n_i = 0$ on $\partial\omega_i$ and
$\psi_{i,j}  \cdot n_{E_i} = \delta^{h}_{j}$ on $E_i$, where $n_{E_i}$ is a unit normal vector on $E_i$,
$\delta^{h}_{j}$ is a fine-scale discrete delta function defined by
\begin{equation}
\delta^{h}_{j} =\{1 \ \mbox{on } e_{j}; \ 0 \  \mbox{on } E_{i}\backslash e_{j}\}
%\begin{cases}
%1 &\mbox{on } e_{j}, \\
%0 &\mbox{on } \partial \omega_{i}\backslash e_{j},
%\end{cases}
%\label{eq:delta_Coupling1}
\end{equation}
and $n$ is the outward unit-normal vector on $\partial \omega_i$. The function $\alpha_{i,j}$ is constant on each coarse grid block within $\omega_i$
and it should satisfy the condition
%$\int_{\omega_i} \alpha_{i,j} = \int_{\partial \omega_{i}} \psi_{i,j}  \cdot n_{i}$.
$\int_{K} \alpha_{i,j} = \int_{\partial K} \psi_{i,j}  \cdot n_{i}$
for each coarse grid $K\in \omega_i$.
We denote the number of snapshots from (\ref{localEquation_coupling})
by $L_i$ and the space spanned by these functions by $V_{H,\text{snap}}^{\omega_i}$.

Now, we define the local spectral problem for $V_{H,\text{snap}}^{\omega_i}$
for the construction of the offline space $V_{H,\text{off}}^{\omega_i}$.
The local spectral problem is to find real number $\lambda$
and $v \in V_{H,\text{snap}}^{\omega_i}$ such that
\begin{equation}\label{spectralProblem}
a_i(v, w) = \lambda s_i(v, w), \qquad \forall w \in V_{H,\text{snap}}^{\omega_i},
\end{equation}
where $a_i$ is a symmetric non-negative definite bilinear operator
and $s_i$ is a symmetric positive definite bilinear operators defined on $V_{H,\text{snap}}^{\omega_i} \times V_{H,\text{snap}}^{\omega_i}$.
Motivated by the analysis (\cite{chung2015mixed}), we choose
\begin{equation*}
a_i(v, w) = \int_{E_i} (v\cdot n_{E_i}) (w\cdot n_{E_i}), \ \
s_i(v, w) = \int_{\omega_{i}} \kappa^{-1} v \cdot w +\text{div}(v)\text{div}( w).
\end{equation*}
We let $\lambda_j^{\omega_i}$ be the eigenvalues of (\ref{spectralProblem}) arranged in ascending order,
and $\psi_j^{\omega_i}$ be the corresponding eigenfunctions.
We will use the first $l_i$ eigenfunctions to construct our offline space $V_{H,\text{off}}^{\omega_i}$.
We note that it is important to keep the eigenfunctions with small eigenvalues in the offline space.
The global offline space $V_{H,\text{off}}$ is
the union of all $V_{H,\text{off}}^{\omega_i}$.

The {\it oversampling} can be performed by using the snapshots
in the regions surrounding the edge and computing boundary
conditions. We refer to
\cite{chung2015mixed} for the details.

\subsubsection{Multiscale basis functions in  GMsDGM. Oversampling}
\label{sec:mixed_ipdg}

For each $K_i\in \mathcal{T}^H$, we consider an oversampled region $K^+_i \supset K_i$. We define $\psi^+_{i,j}$ such that
\begin{equation}
\label{GMsDGM_Nonover_snap}
\begin{aligned}
-\nabla \cdot ( \kappa \nabla\psi^+_{i,j})  =0\ \text{ in } \ K^+_{i},\quad
\psi_{i,j}  =\delta_{j}^{h} \ \text{ on }\ \partial K^+_{i},
\end{aligned}
\end{equation}
where $\delta^{h}_j$ is the delta function defined in Section \ref{sec:gmsfem_basics}.
%is piecewise
%linear on $\partial K^+_i$ with respect to the fine grid such that $\delta^{h}_j$ has the value one at the j-th fine grid node and
%value zero at all the remaining fine grid nodes.
We denote the number of local snapshots by $L^+_i$ and the space spanned by these functions by $V^{+,K_i}_{H,\text{snap}}$.

These basis functions $\psi^+_{i,j}$ are supported in $K^+_i$. We denote the restriction of $\psi^+_{i,j}$ on $K_i$ by   $\psi_{i,j}$. Then we remove the linear dependence of $\psi_{i,j}$ by performing Proper Orthogonal Decomposition (POD)
(\cite{volkwein05}).
Next, we will define the local spectral problem as finding
$\lambda$ and $v\in V^{+,\omega_i}_{H,\text{snap}}$ such that
\begin{equation}
\label{eq:spectral_DG_over}
\begin{aligned}
 a_i(v,w) = \lambda s_i(v,w)\quad  \forall w\in V^{+,\omega_i}_{H,\text{snap}},
\end{aligned}
\end{equation}
where $a_i$ is a symmetric non-negative definite bilinear operator
and $s_i$ is a symmetric positive bilinear operators defined on $V^{+,\omega_i}_{H,\text{snap}} \times V^{+,\omega_i}_{H,\text{snap}}$,
%Motivated by the analysis, we choose
%\begin{equation*}
%\begin{aligned}
$a_i(v, w) = \int_{K_i^+} \kappa \nabla v\cdot \nabla w$ and
$s_i(v, w) = \int_{\partial K_i} \kappa \,v \,w$.

We let $\lambda_j^{\omega_i}$ be the eigenvalues of (\ref{eq:spectral_DG_over}) arranged in ascending order,
and $\psi_j^{+,\omega_i}$ be the corresponding eigenfunctions. We denote $\psi_j^{\omega_i}$ as the restriction of $\psi_j^{+,\omega_i}$ on $K_i$
We will use the first $l_i$ restricted eigenfunctions to construct
the offline space $V_{H,\text{off}}^{\omega_i}$.
We note that it is important to keep the eigenfunctions
with small eigenvalues in the offline space.
The global offline space $V_{H,\text{off}}$ is
the union of all $V_{H,\text{off}}^{\omega_i}$.

\begin{remark}[Energy minimizing oversampling]
\label{rem:energy_min}
For efficient online simulations, one needs to develop
energy minimizing snapshot vectors \cite{chung2015_EMO}.
In these approaches, the snapshot vectors are computed by
solving local constrained minimization problems. More precisely,
first, we use an oversampled region
 to construct the snapshot space. We denote them
$\psi_1^+,...,\psi_N^+$ for simplicity and consider one coarse block.
 Next, we consider the restriction of these snapshot vectors in a target coarse block and identify linearly independent components, $\psi_1,...,\psi_M$, $M\leq N$. A next important ingredient is to construct minimum energy snapshot vectors that represent the snapshot functions in the target coarse region. These snapshot vectors are constructed by solving a local minimization problem.
We denote them
by $\widetilde{\psi}_1,...,\widetilde{\psi}_M$. In the final step, we
perform a local spectral decomposition to compute multiscale basis functions
%using an inner product matrix $A$, which comes
%from the local fine-grid discretization
(\cite{chung2015_EMO}).
Energy minimizing snapshots are important for online basis computations
and are used to show the stable decomposition.
%In particular, we use
%$\frac{(\widetilde{\psi^+})^T A \widetilde{\psi}^+}{\psi^T A \psi}$.
%This general concept can be used in various discretizations and we
%discuss a mixed formulation and Interior Penalty Discontinuous Galerkin (IPDG).
\end{remark}

\subsection{Numerical results}

We present the numerical results for the mixed GMsFEM. The permeability
field is as shown in Figure \ref{fig:perm_random} and
$f=1$. The coarse grid $10\times 10$ and we use zero Dirichlet boundary
conditions. In Table \ref{tab:GMsFEM_mixed}, the relative errors measured in the norm $\|v\|^2=\int_{\Omega} \kappa^{-1} |v|^2$
are shown as we increase the number of basis functions. As we observe from
this table,  the mixed GMsFEM has a very good convergence rate.
In general, we have observed a good accuracy when the
mixed formulation of the GMsFEM is used.
With only $3$ basis functions per node, the error in the velocity field
is less than
$5$\%. We note that the total degrees of freedom for the fine-grid
solution is $20200$.

% \newpage
% \begin{table}[h]
% \centering
% \begin{tabular}{|c|c|c|}
% \hline
% \#basis & ${\|v_{h}-v_{H}\|}/{\|v_{h}\|}$ \tabularnewline
% \hline
% 1 (220) & 16.50\% \%\tabularnewline
% \hline
% 2 (440) & 6.41\% \%\tabularnewline
% \hline
% 3 (660) & 3.91\%\%\tabularnewline
% \hline
% 4 (880) & 2.65\% \%\tabularnewline
% \hline
% 5 (1100) & 1.64\%\%\tabularnewline
% \hline
% \end{tabular}
% \protect\caption{Convergence history with increasing number of offline basis functions for the mixed GMsFEM.}
% \label{tab:GMsFEM_mixed}
% \end{table}

\begin{table}[h]
\centering
\begin{tabular}{|c|c|c|c|c|c|}
\hline
\#basis & 1 (220) & 2 (440) & 3 (660) & 4 (880) & 5 (1100) \tabularnewline
\hline
${\|v_{h}-v_{H}\|}/{\|v_{h}\|}$ &
16.50\% &  6.41\%  & 3.91\% & 2.65\% &  1.64\%  \tabularnewline
\hline
\end{tabular}
\protect\caption{Convergence history with increasing number of offline basis functions for the mixed GMsFEM.}
\label{tab:GMsFEM_mixed}
\end{table}

% \begin{table}[h]
% \centering
% \begin{tabular}{|c|c|c|}
% \hline
% \#basis & $\cfrac{\|v_{h}-v_{ms}\|}{\|v_{h}\|}$ & $\cfrac{\|p_{h}-p_{ms}\|}{\|p_{h}\|}$\tabularnewline
% \hline
% 1 (220) & 16.50\% & 14.97\%\tabularnewline
% \hline
% 2 (440) & 6.41\% & 14.66\%\tabularnewline
% \hline
% 3 (660) & 3.91\% & 14.65\%\tabularnewline
% \hline
% 4 (880) & 2.65\% & 14.65\%\tabularnewline
% \hline
% 5 (1100) & 1.64\% & 14.65\%\tabularnewline
% \hline
% \end{tabular}
% \protect\caption{Convergence history with increasing number of offline basis functions, fine grid 20200}
% \label{tab:GMsFEM_mixed}
% \end{table}

\section{Discussions on the sparsity in the snapshot space}
\label{sec:sparsity}

In previous approaches,
multiscale basis functions are constructed using local snapshot spaces,
where a snapshot space is a large space that represents the solution
behavior in a coarse block. In a number of applications,
one may have a sparsity in the snapshot space for an appropriate choice
of a snapshot space. More precisely,
the solution may only involve a portion of the snapshot space.
In this case, one can use sparsity techniques
to identify multiscale basis
functions. We briefly discuss two such sparse local
multiscale model reduction approaches (see \cite{chung2015sparse} for details).

In the first approach (which is used for parameter-dependent multiscale PDEs),
we use local minimization techniques, such as sparse POD, to identify multiscale
basis functions, which are sparse in the snapshot space.
These minimization techniques use $l_1$ minimization to find
local multiscale basis functions, which are further used for finding
the solution. In the second approach (which is used
for the Helmholtz equation), we directly apply $l_1$
minimization techniques to solve the underlying PDEs. This approach is more
expensive as it involves a large snapshot space; however, in this example,
we cannot identify a local minimization principle, such as local generalized
SVD.

We consider
$\mathcal{L} u = f$,
where $\mathcal{L}$ is a differential operator. For example,  we consider
parameter-dependent heterogeneous flows, $\mathcal{L}u=
-\text{div}(\kappa(x;\mu)\nabla u)$,
 and the Helmholtz equation, $\mathcal{L}u=-\text{div}(\kappa(x)\nabla u) - \zeta^2 n(x) u$.
Previous approaches
use all snapshot vectors when seeking multiscale basis functions. In a number of
applications, the solution is sparse in the snapshot space, which implies
that in the expansion
\[
u=\sum_{i,j} c_{i,j} \psi_i^{\omega_j},
\]
many coefficients $c_{i,j}$ are zeros.
In this case, one can save computational effort
by employing sparsity techniques.

\begin{figure}[htb]
  \centering
  \includegraphics[width=0.65 \textwidth]{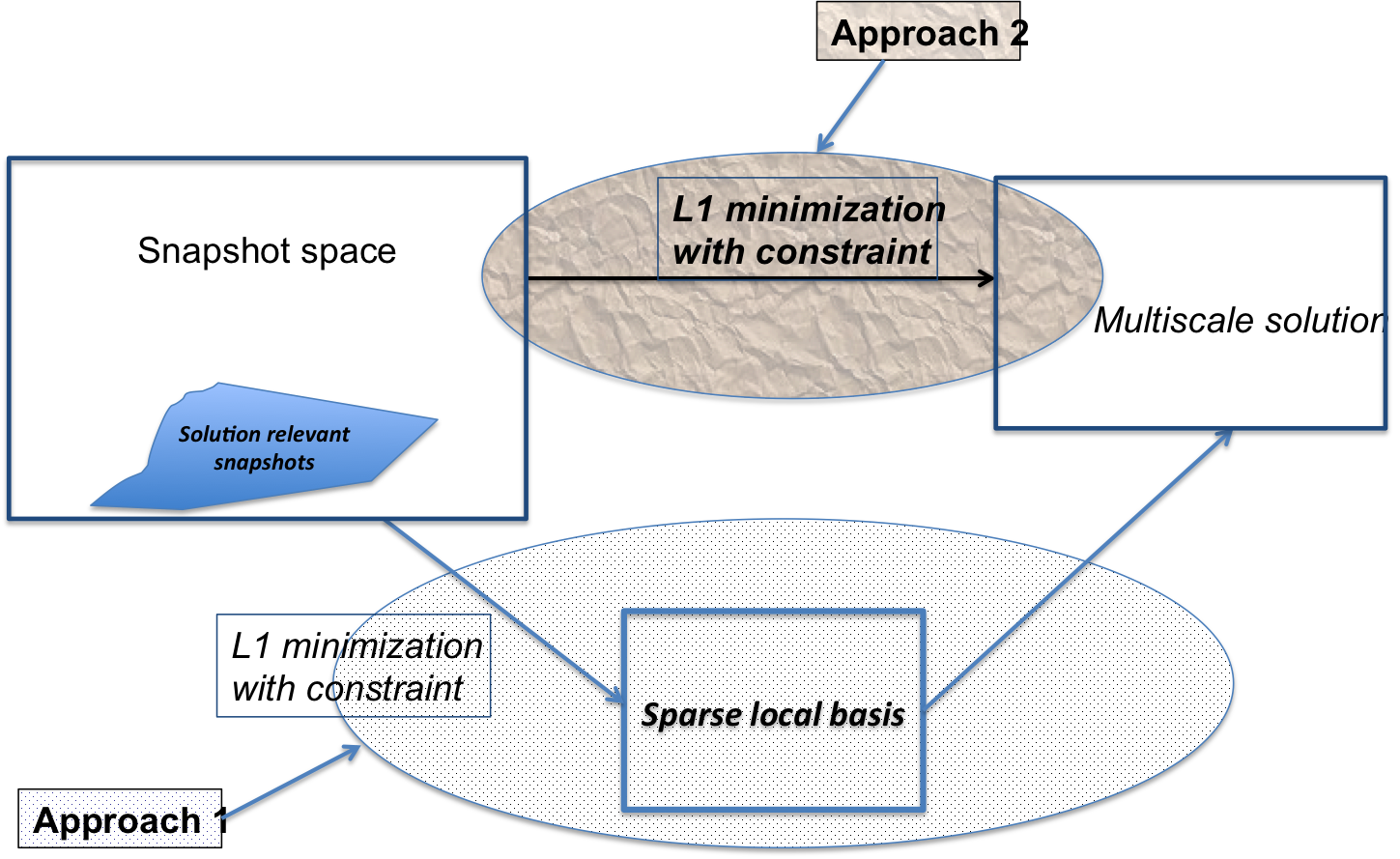}
  \caption{Illustration of the approaches.}
  \label{overview1}
\end{figure}

The main challenge in these applications is to construct a snapshot space,
where the solution is sparse.
In the first example, this can be achieved,
because an online parameter value $\mu$ can be close to some of the pre-selected
 offline
values of $\mu$'s, and thus, the multiscale basis functions (and the solution)
can have a sparse representation in the snapshot space. In the
second example,
we select cases where the solution $u$ contains only a few snapshot
vectors corresponding to some wave directions.
We note that if the snapshot space is not chosen carefully, one may not have
the sparsity.
We can consider two distinct cases.
\begin{itemize}

\item First approach: ``Local-Sparse Snapshot Subspace Approach''. Determining the online sparse space locally via local spectral sparse decomposition in the snapshot space (motivated by
parameter-dependent problems).

\item  Second approach: ``Sparse Snapshot Subspace Approach''.  Determining the online space globally via a global solve  (motivated by using plane wave snapshot vectors and the Helmholtz equation).

\end{itemize}
See Figure \ref{overview1} for illustration.
We use sparsity techniques
(e.g., \cite{candes2006compressive,schaeffer2013sparse})
to identify local multiscale basis functions
and solve the global problem. The numerical results and more discussions
can be found in \cite{chung2015sparse}.

\section{Space-time GMsFEM}
\label{sec:spacetime}

%In this section, we discuss the GMsFEM for problems with space-time
%heterogeneities.
Many multiscale processes vary over multiple space and time scales.
In many applications,
the heterogeneities change due to the time can be
significant and it needs to be taken into account in reduced-order
models.
Some well-known approaches for handling {\it separable}
spatial and temporal scales
are homogenization techniques
\cite{jikov2012homogenization,pankov2013g,pavliotis2008multiscale,efendiev2005homogenization}.
 In these methods, one solves local problems
in space and time. We remind
a well-known case of the parabolic equation
\begin{equation}
\label{eq:para PDE}
\begin{split}
{\partial u/\partial t}
-\text{div}(\kappa(x,t)\nabla u)  = f,
\end{split}
\end{equation}
subject to smooth initial and boundary conditions.
The homogenized equation has the same form
as (\ref{eq:para PDE}), but with the smooth coefficients
$\kappa^*(x,t)$. One can compute the coefficients using the solutions of local space-time parabolic equations in the  periodic cell.
This localization is possible thanks to the scale separation.
One can extend this homogenization procedure to numerical homogenization
type methods
\cite{ming2007analysis,fish_book}, where
one solves the local parabolic equations in each coarse block.

Previous approaches
within GMsFEM focused
on constructing
 multiscale spaces and relevant ingredients in space only.
The extension of the GMsFEM to space-time heterogeneous problems
requires a modification of the space only problems because
of (1) the parabolic nature of cell solutions, (2) extra degrees of freedom
associated with space-time cells,
and (3) local boundary conditions in space-time cells.
In the approach that we are going to discuss,
we construct snapshot spaces in space-time
local domains. We construct the snapshot solutions
by solving local problems.
We can construct a complete snapshot space by taking all
possible boundary conditions; however, this can lead to very high
computational cost. For this reason, we
use randomized boundary conditions
for local snapshot vectors.
 Computing multiscale basis functions employs local
spectral problems in space-time domain.

%\subsection{Space-time GMsFEM}\label{sect:GMsFEM}

We consider the  parabolic differential equation
(\ref{eq:para PDE}) in a space-time domain  $\Omega \times (0,T)$
%and the time
%interval $[0,T]$
and assume $u=0$ on $\partial\Omega\times(0,T)$ and
$u(x,0)=\beta(x)$ in $\Omega$.
The proposed method follows the space-time finite element framework, where the time dependent multiscale basis functions are constructed on the coarse grid. Therefore, compared with the time independent basis structure, it gives a more efficient numerical solver for the parabolic problem in complicated media.

We would like to compute the solution $u_{H}$ in the whole
time interval $(0,T)$. In fact, if we assume the solution
space $V_{H,\text{off}}^{(0,T)}$ is a direct sum of the spaces only containing the functions
defined on one single coarse
time interval $(T_{n-1},T_{n})$, we can decompose
the problem into a sequence of problems
and find the solution $u_{H}$ in each time interval sequentially.
The coarse space will be constructed in each time interval and
$V_{H,\text{off}}^{(0,T)}=\oplus_{n=1}^{N}V_{H,\text{off}}^{(T_{n-1},T_n)}$,
 where $V_{H,\text{off}}^{(T_{n-1},T_n)}$ only contains the functions having zero values
in the time interval $(0,T)$ except $(T_{n-1},T_{n})$, namely
$\forall v\in V_{H,\text{off}}^{(T_{n-1},T_n)},$
$v(\cdot,t)=0\text{ for }t\in(0,T)\backslash(T_{n-1},T_{n})$.

The coarse-grid equation
consists of finding
 $u_{H}^{(n)}\in V_{H,\text{off}}^{(T_{n-1},T_n)}$
(where $V_{H,\text{off}}^{(T_{n-1},T_n)}$ will be defined later)
satisfying
\begin{equation}
\label{eq:space-time-FEM-coarse-decoupled}
  \int_{T_{n-1}}^{T_{n}}\int_{\Omega}\cfrac{\partial u_{H}^{(n)}}{\partial t}v+\int_{T_{n-1}}^{T_{n}}\int_{\Omega}\kappa\nabla u_{H}^{(n)}\cdot\nabla v+\int_{\Omega}u_{H}^{(n)}(x,T_{n-1}^{+})v(x,T_{n-1}^{+})\nonumber
=  \int_{T_{n-1}}^{T_{n}}\int_{\Omega}fv+\int_{\Omega}g_{H}^{(n)}(x)v(x,T_{n-1}^{+}),
\end{equation}
for all $v\in V_{H,\text{off}}^{(T_{n-1},T_n)}$,
where
%\[
$g_{H}^{(n)}(\cdot)=\{u_{H}^{(n-1)}(\cdot,T_{n-1}^{-}) \text{ for }n\geq1;
\beta(\cdot)  \text{ for }n=0\}$,
and $F(\alpha^+)$ and $F(\alpha^-)$ denote the right hand and left hand limits of $F$ at $\alpha$ respectively.
%g_{H}^{(n)}(\cdot)=\begin{cases}
%u_{H}^{(n-1)}(\cdot,T_{n-1}^{-}) & \text{ for }n\geq1,\\
%\beta(\cdot) & \text{ for }n=0.
%\end{cases}
%\]
Then, the solution $u_{H}$ of the problem in $(0,T)$
is the direct sum of all these $u_{H}^{(n)}$'s, that is $u_{H}=\oplus_{n=1}^{N}u_{H}^{(n)}$.

Next, we motivate the use of space-time multiscale basis functions
by comparing it to space multiscale basis functions.
In particular, we discuss the savings in the reduced models when space-time
multiscale basis functions are used compared to space multiscale
basis functions.
We denote $\{t_{n1},\cdot\cdot\cdot,t_{np}\}$ as the $p$ fine time steps
in $(T_{n-1},T_n)$. When we construct space-time multiscale basis functions,
the solution can be represented as $u_H^{(n)} = \sum_{l,i} c_{l,i} \psi_l^{\omega_i}(x,t)$
in the interval  $(T_{n-1},T_n)$. In this case, the number of coefficients
$c_{l,i}$ is related to the size of the reduced system in space-time interval.
On the other hand, if we use only space multiscale basis functions,
we need to construct these multiscale basis functions at each
fine time instant $t_{nj}$, denoted by $\psi_{l}^{\omega_i}(x,t_{nj})$.
The solution $u_H$ spanned by these basis functions will have a much
larger dimension because each time instant is represented by
multiscale basis functions. Thus, performing space-time multiscale
model reduction can provide a substantial CPU savings.

\subsection{Construction of offline basis functions}
\subsubsection{Snapshot space}

Let $\omega$ be a given coarse neighborhood in space.
We omit the coarse node index to simplify the notations. The construction of the
offline basis functions on coarse time interval $(T_{n-1},T_n)$ starts
with a snapshot space $V_{H,\text{snap}}^{\omega}$ (or $V_{H,\text{snap}}^{\omega,(T_{n-1},T_n)}$). We also omit the coarse time index $n$ to simplify the notations. The snapshot space
$V_{H,\text{snap}}^{\omega}$ is a set of functions defined on $\omega$
and contains all or most necessary components of the fine-scale
solution restricted to $\omega$. A spectral problem is then solved
in the snapshot space to extract the dominant modes in the snapshot space.
These dominant modes are the offline basis functions and the resulting
reduced space is called the offline space. There are two choices of
 $V_{H,\text{snap}}^{\omega}$ that are commonly used.

The first choice is to use all possible fine-grid functions in
$\omega\times (T_{n-1},T_{n})$.
This snapshot space provides an accurate
approximation for the solution space; however,
this snapshot space can be very large.
The second choice for the snapshot spaces consists
of solving local problems for all possible boundary conditions,
which we present here.
As before,
we denote by $\omega^{+}$ the oversampled space region of
$\omega \subset\omega^{+}$, defined by adding several fine- or coarse-grid
layers around $\omega$. Also, we define $(T_{n-1}^{*}, T_{n})$ as
the left-side oversampled time region for $(T_{n-1},T_{n})$. In the following,
we generate inexpensive snapshots using random boundary conditions on
the oversampled space-time region $\omega^{+}\times(T_{n-1}^{*},T_{n})$.
That is for each fine boundary
node on $\partial \left( \omega^+\times (T_{n-1}^*,T_{n}) \right)$, we solve
a small number of local problems imposed with random boundary conditions
\begin{equation*}
\begin{split}
\partial\psi_{j}^{+,\omega}/\partial t -\text{div} (\kappa(x,t) \nabla \psi_{j}^{+,\omega})=0\ \ \text{in}\ \omega^{+}\times (T_{n-1}^{*},T_{n}), \ \
\psi_{j}^{+,\omega}(x,t)= r_l\  \ \text{on} \ \ \partial \left( \omega^{+}\times (T_{n-1}^{*},T_{n}) \right),
\end{split}
\end{equation*}
where $r_l$ are independent identically distributed (i.i.d.) standard Gaussian random vectors on the fine-grid nodes of the boundaries $t=T_{n-1}^{*}$ and on $\partial \omega^{+}\times (T_{n-1}^{*},T_{n})$.

Then the local snapshot space on $\omega^{+}\times (T_{n-1}^{*},T_{n})$ is
\[
V_{H,\text{snap}}^{+,\omega} = \text{span}\{\psi_{j}^{+,\omega}(x,t) | j=1,\cdot\cdot\cdot, l^{\omega}+p_{\text{bf}}^{\omega}\},
\]
where $l^{\omega}$ is the number of local offline basis we want to construct in $\omega$ and $p_{\text{bf}}^{\omega}$ is the buffer number. Later on, we use the same buffer number for all $\omega$'s and simply use the notation $p_{\text{bf}}$. In the following sections, if we specify one special coarse neighborhood $\omega_i$, we use the notation $l_i$ to denote the number of local offline basis. With these snapshots, we follow the procedure in the following subsection to generate offline basis functions by using an auxiliary spectral decomposition.

%\subsubsection{Offline space}
\subsubsection{Offline space}

To obtain the offline basis functions, we need to perform a space reduction by appropriate spectral problems. Motivated by a convergence analysis, we adopt the following spectral problem on $\omega^{+}\times (T_{n-1},T_{n})$.
Find $(\phi,\lambda)\in V_{H,\text{snap}}^{+,\omega}\times\mathbb{R}$ such that
\begin{equation}\label{eq:eig-problem}
A_n(\phi,v) = \lambda S_n(\phi,v), \quad \forall v \in V_{\text{snap}}^{\omega^{+}},
\end{equation}
where the bilinear operators $A_n(\phi,v)$ and $S_n(\phi,v)$ are defined by
\begin{equation}
\begin{split}
A_n(\phi,v) = \frac{1}{2} \left( \int_{\omega^{+}}\phi(x,T_{n})v(x,T_{n}) + \int_{\omega^{+}}\phi(x,T_{n-1})v(x,T_{n-1}) \right) +
 \int_{T_{n-1}}^{T_{n}}\int_{\omega^{+}}\kappa(x,t)\nabla\phi \cdot \nabla v, \\
S_n(\phi,v) = \int_{\omega_{+}}\phi(x,T_{n-1})v(x,T_{n-1}) +\int_{T_{n-1}}^{T_{n}}\int_{\omega^{+}}\widetilde{\kappa}^{+}(x,t)\phi v,
\end{split}
\end{equation}
where the weight function $\widetilde{\kappa}^{+}(x,t)$ is defined by
$\widetilde{\kappa}^{+}(x,t) = \kappa(x,t)\sum_{i=1}^{N_c}|\nabla\chi_i^{+}|^2$,
$\{\chi_i^{+}\}_{i=1}^{N_c}$ is a partition of unity associated with the oversampled coarse neighborhoods $\{\omega_i^{+}\}_{i=1}^{N_c}$ and satisfies $|\nabla\chi_i^{+}|\geq|\nabla\chi_i|$ on $\omega_i$,
where $\chi_i$ is the standard multiscale basis function for the coarse node $x_i$,
% (that is, with linear boundary conditions for cell problems). More precisely,
$-\text{div}(\kappa(x,T_{n-1})\nabla\chi_i) = 0$, in $K$,
 $\chi_i=g_i$ on $\partial K$,
for all $K\in\omega_i$, where $g_i$
%is a continuous function on
%$\partial K$ and
is linear on each edge of $\partial K$.

We arrange the eigenvalues $\{\lambda_j^{\omega}|j=1,2,\cdot\cdot\cdot\,L_{\omega}+p_{\text{bf}}^{\omega}\}$ from (\ref{eq:eig-problem}) in the ascending order, and select the first $L_{\omega}$ eigenfunctions, which are corresponding to the first $L_{\omega}$ ordered eigenvalues, and
then we can obtain the dominant modes $\psi_{j}^{\omega}(x,t)$ on the target
region $\omega\times (T_{n-1},T_{n})$ by restricting
$\psi_j^{+,\omega}(x,t)$ onto $\omega\times (T_{n-1},T_{n})$. Finally, the offline basis functions on $\omega\times (T_{n-1},T_{n})$ are defined by $\phi_j^{\omega}(x,t) = \chi^\omega\psi_j^{\omega}(x,t)$, where $\chi^\omega$ is the standard multiscale basis function  for a generic coarse neighborhood $\omega$. This product gives conforming basis functions (Discontinuous Galerkin discretizations can also be used). We also define the local offline space on $\omega\times (T_{n-1},T_{n})$ as
\[
V_{H,\text{off}}^{\omega} = \text{span}\{\phi_{j}^{\omega}(x,t) | j=1,\cdot\cdot\cdot, l^{\omega} \}.
\]

Note that one can take $V_{H,\text{off}}^{(T_{n-1},T_n)}$ in the coarse-grid
equation
% in
%(\ref{eq:space-time-FEM-coarse-decoupled})
as $V_{H,\text{off}}^{(T_{n-1},T_n)} =  \text{span}\{\phi_{j}^{\omega_i}(x,t) | 1\leq i\leq N_c, 1\leq j\leq l_{i} \}$.
% As a result, $V_{H} = V_{\text{off}}= \oplus_{n=1}^{N}V_{H}^{(n)}$.

%%%%%STOPPED HERE%%%%%%%%%%%%%%%%%

\subsection{Numerical result}

We start with a high-contrast permeability field shown in Figure \ref{fig:perm_random}, which is translated uniformly in $x_2$ direction
 after every other fine time step. The total translation is $10$\% of the global
domain.    The number of local offline basis that will be used in each $\omega_i$, is denoted by $l_i$, and the buffer number $p_{\text{bf}}$ needs to be chosen in advance since they determine how many local snapshots are used. Then, we can construct the lower dimensional offline space by performing space reduction on the snapshot space. In the experiments, we use the same  buffer number and the same number of local offline basis for all coarse neighborhood $\omega_i$'s.
Our numerical experiments show that the error does not change
beyond
$p_{\text{bf}}=4$. In our numerical simulations, we take $p_{\text{bf}}=8$
and vary $l_i$. The errors are displayed in  Table \ref{Table:test bf Li 1}.
These are $L^2$ and energy errors.  We observe that with a fixed buffer number, the relative errors are decreasing as using more offline basis.
We can also observe that the values
$1/\Lambda_{*}$'s and energy errors are
 correlated as predict our theory \cite{chung2015_spacetime}.

% \begin{table}[h]
% \centering
%  \begin{minipage}[t]{0.4\linewidth}
%  \centering
%   \begin{tabular}{ | c | c | c | c |  }
%     \hline
%     $p_{\text{bf}}$  &Snapshot ratio  &$e_1$  &$e_2$  \\
%     \hline
%     1   & 0.0158     & 6.18\%    & 53.90\%   \\
%     4   & 0.0197     & 5.66\%    & 48.04\%   \\
%     8   & 0.0250     & 5.17\%    & 45.86\%   \\
%     12  & 0.0302     & 5.16\%    & 43.83\%   \\
%     20  & 0.0407     & 4.71\%    & 41.14\%   \\
%     30  & 0.0539     & 4.35\%    & 38.68\%   \\
%     40  & 0.0670     & 4.23\%    & 37.60\%   \\
%     \hline
%   \end{tabular}
%   \end{minipage}
%   \begin{minipage}[t]{0.5\linewidth}
%   %\centering
%   \begin{tabular}{ | c | c | c | c | c |  }
%     \hline
%     $L_i$  &$dim(V_{\text{off}})$   &Snapshot ratio  &$e_1$  &$e_2$  \\
%     \hline
%     2     &162   & 0.0131   &17.03\%    & 129.14\%   \\
%     6     &486   & 0.0184   &8.11\%    & 62.59\%  \\
%     10    &810   & 0.0237   &6.97\%    & 54.85\%   \\
%     20    &1620  & 0.0368   &4.81\%    & 41.18\%   \\
%     30    &2430  & 0.0499   &3.29\%    & 31.64\%   \\
%     40    &3240  & 0.0631   &2.28\%    & 24.43\%    \\
%     50    &4050  & 0.0762   &1.54\%    & 18.45\%   \\
%     \hline
%   \end{tabular}
%   \end{minipage}
% \caption{First permeability field. Left: errors with the fixed number of offline basis $L_i=11$. Right: errors with the fixed buffer number $p_{\text{bf}}=8$. }
% \label{Table:test bf Li 1}
% \end{table}

\begin{table}[h]
\centering
  %\centering
  \begin{tabular}{ | c | c | c | c | c |  }
    \hline
    $L_i$  &$dim(V_{\text{off}})$   &$L^2([0,T],\Omega)$  & energy  \\
    \hline
    2     &162      &23.05\%    & 87.86\%   \\
    6     &486      &7.86\%    & 43.29\%  \\
    10    &810      &5.93\%    & 35.23\%   \\
    14    &1134     &3.43\%    & 25.66\%   \\
    18    &1458     &1.57\%    & 17.42\%   \\
    22    &1782     &0.88\%    & 12.89\%    \\
    \hline
  \end{tabular}
\caption{Relative errors using different $l_i$'s. }
\label{Table:test bf Li 1}
\end{table}

\section{GMsFEM in perforated domains}
\label{sec:perforated}

One important class of multiscale problems consists
in perforated
domains. In these problems, differential equations are
formulated in perforated domains. These domains
can be considered to be the outside of inclusions or connected
bodies of various sizes.
Due to the variable sizes and geometries of these perforations,
solutions to these problems have multiscale features.
One solution approach
involves posing the problem in a domain without perforations
but with a very high contrast penalty term representing the domain
heterogeneities
(\cite{griebel2010homogenization}).
However, the void space can be a small portion of the
whole domain and, thus, it is computationally expensive to
enlarge the domain substantially.

 The main difference in developing
multiscale methods for problems in perforated domains is the complexity
of the domains and that many portions of the domain are excluded
in the computational domain. This poses a challenging task.
For typical
upscaling and numerical homogenization
(e.g., \cite{sanchez1980non, henning2009heterogeneous}),
the macroscopic equations
do not contain perforations and one computes the effective properties.

Several multiscale methods have been developed for
problems in perforated domains. The use of GMsFEM for solving
multiscale problems in perforated domains is motivated
by recent works \cite{Maz'ya13, Bris14, henning2009heterogeneous}.
Using snapshot spaces
is essential in problems with perforations, because the snapshots contain
necessary geometry information. In the snapshot space, we perform
local spectral decomposition to identify multiscale basis functions.
In this section, we show an example for Stokes
equations and refer to \cite{chung2015generalizedperforated}
for more discussions and extensions to other problems.
%We also present some results for the online computations
%and the numerical results.

%structure
%\subsection{Preliminaries}
%\label{prelim}
%==============================

\subsection{Problem setting}

In this section, we present the underlying problem as stated in \cite{CELV2015}
and the corresponding fine-scale and coarse-scale discretization.
Let $\Omega\subset \mathbb{R}^d$ ($d=2,3$) be a bounded domain covered
by inactive cells (for Stokes flow and Darcy flow) or active cells (for elasticity problem) $\mathcal {B}^{\eps}$.
We will consider $d=2$ case, though the results can be extended
to $d>2$.
We use the superscript $\eps$ to denote quantities related to perforated domains.
The active cells are where the underlying problem is solved, while inactive cells are the rest of the region.
Suppose the distance between inactive cells (or active cells) is
of order $\eps$. Define $\Omega^{\eps}:=\Omega\backslash \mathcal {B}^{\eps}$. See Figure \ref{fig:twodomain} (left)
for an illustration of the perforated domain.
We consider the following problem defined in a perforated domain $\Omega^{\eps}$
\begin{equation}
\label{eq:original_perforated}
\mathcal {L}^{\eps}(w)=f \  \text{in} \ \ \Omega^{\eps},\ \
w=0 \text{ or } \frac{\partial w}{\partial n}=0 \text{ on }\partial \Omega^{\eps}\cap  \partial \mathcal {B}^{\eps},
\end{equation}
$w=g$ on $\partial \Omega\cap \partial \Omega^{\eps}$,
where $\mathcal {L}^{\eps}$ denotes a linear differential operator, $n$ is the unit outward normal to the boundary, $f$ and $g$ denote given functions with sufficient regularity.
We will focus on the Dirichlet problem, namely, $w=0$ in (\ref{eq:original_perforated}).

%\begin{figure}[htb]
%  \centering
%  \includegraphics[width=0.3 \textwidth]{domain}
%  \caption{Illustration of a perforated domain.}
%  \label{fig:perf_domain}
%\end{figure}

Denote by ${V}(\om^{\eps})$ the appropriate solution space, and
 ${V}_{0}(\om^{\eps})=\{v\in {V}(\om^{\eps}), v=0 \text{ on }\partial\om^{\eps}\}$.
The variational formulation of \eqref{eq:original_perforated} is to find $w\in {V}(\om^{\eps})$ such that
\[
\avrg{\mathcal {L}^{\eps}(w),v}{\om^{\eps}}= (f,v)_{\om^{\eps}} \qquad  \text{for all } v \in V_{0}(\om^{\eps}),
\]
where $\avrg{\cdot,\cdot}{\om^{\eps}}$ denotes a specific inner product over $\om^{\eps}$ for either scalar functions or vector functions and
and $(f,v)_{\om^{\eps}}$ is the $L^2$ inner product.
Some specific examples for the above abstract notations are given below.

\textbf{Laplace:} For the Laplace operator with homogeneous Dirichlet boundary conditions on $\partial \om^{\eps}$, we have
$\mathcal {L}^{\eps}(u)=-\Delta u$,
and ${V}(\om^{\eps})=H^{1}_{0}(\Omega^{\eps})$, $\avrg{\mathcal {L}^{\eps}(u),v}{\om^{\eps}}=(\nabla u,\nabla v)_{\om^{\eps}}$.

\textbf{Elasticity:} For the elasticity operator with a
homogeneous Dirichlet boundary condition on $\partial \om^{\eps}$, we assume the medium is isotropic.
Let $ {u}\in (H^{1}(\Omega^{\eps}))^{2}$ be the displacement field.
The strain tensor $ {\strain}( {u})\in (L^{2}(\Omega^{\eps}))^{2\times 2}$
is defined by
$ {\strain }( {u}) = \frac{1}{2} ( \nabla  {u} + \nabla  {u}^T )$.
Thus, the stress tensor $ {\sigma}( {u})\in (L^{2}(\Omega^{\eps}))^{2\times 2}$ relates to the strain tensor $ {\strain}( {u})$ such that
$ {\sigma}(u)= 2\mu  {\strain} + \xi \nabla\cdot  {u} \,  {I}$,
where $\xi >0$ and $\mu>0$ are the Lam\'e coefficients. We have
%\begin{align}\label{eqn:elasticity}
$\mathcal {L}^{\eps}(u)=- \nabla \cdot  {\sigma}$,
%\end{align}
where ${V}(\om^{\eps})=(H^{1}_{0}(\Omega^{\eps}))^{2}$ and
$\avrg{\mathcal {L}^{\eps}(u),v}{\om^{\eps}}
=2\mu(\strain(u),\strain(v))_{\om^{\eps}}+\xi (\nabla\cdot u,\nabla\cdot v)_{\om^{\eps}}$.

\textbf{Stokes:} For Stokes equations, we have
% \begin{align}\label{eqn:Stokes}
% \mathcal {L}^{\eps}(u\;,p)=\begin{pmatrix}
% \nabla p -\Delta {u}\\
% \nabla \cdot {u}
% \end{pmatrix},
% \end{align}
\begin{equation}\label{eqn:Stokes}
\mathcal {L}^{\eps}(u\;,p)= (\nabla p -\Delta {u}, \nabla \cdot {u})^T,
\end{equation}
where $\mu$ is the viscosity, $p$ is the fluid pressure, $u$ represents the velocity, ${V}(\om^{\eps})=(H^{1}_{0}(\Omega^{\eps}))^{2}\times L^{2}_{0}(\Omega^{\eps})$, and
\[\avrg{\mathcal {L}^{\eps}(u\;,p),(v\;,q)}{\om^{\eps}}=
\begin{pmatrix}
(\nabla u,\nabla v)_{\om^{\eps}} &-(\nabla \cdot v,p)_{\om^{\eps}}\\
(\nabla \cdot u,q)_{\om^{\eps}}&0
\end{pmatrix}.\]
We recall that $L^{2}_{0}(\Omega^{\eps})$ contains functions in $L^{2}(\Omega^{\eps})$
with zero average in $\Omega^{\eps}$.
We will illustrate our ideas and results using the Stokes equations.

%\subsubsection{Coarse and fine grid notations}

For the numerical approximation of the above problems, we first introduce the notations of fine and coarse grids pertinent to perforated domains.
It follows similar constructions as before with the exception
of domain geometries, which can intersect the boundaries
of coarse regions.
Let $\mathcal{T}_H$ be a coarse-grid partition of the domain
$\Omega^{\eps}$ with mesh size $H$.
Notice that, the edges of the coarse elements do not necessarily have straight edges
because of the perforations (see Figure \ref{fig:twodomain}, left).
By conducting a conforming
refinement of the coarse mesh $\mathcal{T}^H$, we can obtain
 a fine mesh $\mathcal{T}^h$ of $\Omega^{\eps}$ with mesh size $h$.
 %where we assume that the fine mesh resolves the multiscale perforations.
Typically, we assume that $0 < h \ll H < 1$, and that the fine-scale mesh $\mathcal{T}_h$
is sufficiently fine to fully resolve the small-scale information of the domain, and $\mathcal{T}_H$ is a coarse mesh containing many fine-scale features.
Let $N_v$ and $N_e$ be the number of nodes and edges in coarse grid respectively.
We denote by $\{x_i|1\leq i \leq N_c\}$ the set of coarse nodes, and $\{E_j| 1\leq j\leq N_e \}$ the set of coarse edges.
We define a coarse neighborhood $\omega_i^{\eps}$ for each coarse node $x_i$ by
$\omega_i^{\eps} = \cup\{K_j^{\eps}\in \mathcal{T}^H; ~~x_i\in \bar{K_j^{\eps}} \}$,
which is the union of all coarse elements having the node $x_i$.
For the Stokes problem, additionally, we define a coarse neighborhood $\omega_m^{\eps}$ for each coarse edge $E_m$ by
\begin{equation}
\label{eq:n2}
\omega_m^{\eps} = \cup\{K_j^{\eps}\in \mathcal{T}^H; ~~E_m\in \bar{K_j^{\eps}} \},
\end{equation}
which is the union of all coarse elements having the edge $E_m$.
%See Figure \ref{fig：nbhd_def_perfor} for an illustration of the coarse neighborhoods.
%On the triangulation $\mathcal{T}_h$,
%we introduce the following finite element spaces
%\begin{align*}
%{V}_h  &:= \{ {v} \in V(\Omega^{\eps})| {v}|_K \in (P^k(K))^l
%\mbox{ for all } K \in \mathcal{T}_h \},
%\end{align*}
%where, $P^k$ denotes the polynomial of degree $k$( $k=0,\;1,\;2$), and $l$( $l=1,\;2$)
%indicates either a scalar or a vector.
%Note that for the Laplace and elasticity operators, we choose $k=1$, i.e., piecewise linear function space as the fine-scale approximation space; for Stokes problem, we use $(P^2(K))^2$ for fine-scale velocity approximation and $P^0(K)$ for fine-scale pressure approximation.
%We use $Q_h$ to denote the space for pressure.
We let $V_h$ be the fine scale space for velocity $u$ and $Q_h$ be the fine scale space for the pressure $p$.
One can choose $V_h$ to be piecewise quadratic and $Q_h$ to be piecewise constant with respect to the fine mesh $\mathcal{T}^h$.

We will then obtain the fine-scale solution $(u,p)\in V_h \times Q_h$ for the Stokes system
by solving the following variational problem
\begin{align}\label{eq:fine_system1}
\avrg{\mathcal {L}^{\eps}(u,p),(v,q)}{\om^{\eps}}= ((f,0),(v,q))_{\om^{\eps}}, \qquad  \text{for all } (v,q) \in V_h \times Q_h.
\end{align}
These solutions are used as reference solutions to test the performance of the schemes.

%------------------------------------
\subsection{The construction of offline basis functions}
\label{Construction}
%------------------------------------

In this section, we describe the construction of offline basis for the Stokes problem in perforated domains. We refer to \cite{chung2015_perforated_online} for online basis
construction and the analysis, as well as GMsFEM for the elliptic equation and the elasticity equations in perforated domains.

\subsubsection{Snapshot space}

%A snapshot space is a space which contains an extensive set of basis functions that are solutions
%of local problems with all possible boundary conditions up to fine-grid resolution.
To compute snapshot
functions, we solve the following problem on the coarse neighborhood $\omega_i^{\eps}$:
find $(u^i_l, p^i_l)$ (on a fine grid) such that
\begin{equation}
\label{eq:hrmext}
\begin{split}
\int_{\omega_i^{\eps}}   \nabla u^i_l  : \nabla v
- \int_{\omega_i^{\eps}} p^i_l \text{div}(v)  = 0,\  \forall v \in V^i_{h,0},\ \
\int_{\omega_i^{\eps}} q \text{div}(u^i_l) = \int_{\omega_i^{\eps}} c  q, \ \forall q \in Q^i_h,
\end{split}
\end{equation}
with boundary conditions
$u^i_l  =  (0, 0), \text{ on } \partial \mathcal {B}^{\eps}$,
$u^i_l  = (\delta^h_l, 0)  \text{ or } (0, \delta^h_l) ,  \text{ on } \partial \omega_i^{\eps} \backslash  \partial \mathcal {B}^{\eps}$.
Here, we write $\omega_i^{\eps} \backslash  \partial \mathcal {B}^{\eps} = \cup_{l=1}^{S_i} e_l$,
where $e_l$ are the fine-grid edges and $S_i$ is the number of these fine
grid edges on $\omega_i^{\eps}\backslash \partial \mathcal {B}^{\eps}$.
Moreover, $\delta^h_l$ is a fine-scale delta function such that  it has value
$1$ on $e_l$ and value $0$ on other fine-grid edges.
In (\ref{eq:hrmext}), we define $V_h^i$ and $Q_h^i$ as the restrictions of the fine grid space in $\omega_i^{\eps}$
and $V_{h,0}^i \subset V_h^i$ are functions that vanish on $\partial\omega_i^{\eps}$.
We remark that we choose the constant $c$ in \eqref{eq:hrmext} by compatibility condition, $c = \frac{1}{|\omega_i^{\eps}|} \int_{\partial \omega_i^{\eps} /  \partial \mathcal {B}^{\eps} } u^i_l \cdot n_i \, ds$.
We emphasize that, for the Stokes problem, we will solve (\ref{eq:hrmext})
in both node-based coarse neighborhoods and edge-based coarse neighborhoods (\ref{eq:n2}).
Collecting  the solutions of the local problems generates the snapshot space,
$\psi_l^{\omega_i^\epsilon} = u^i_l$ in $\omega_i^{\eps}$:
\[
V_{H,\text{snap}}^{\omega_i^\epsilon} = \{ \psi_l^{\omega_i^\epsilon}: 1 \leq l \leq 2 S_i, \, 1 \leq i \leq (N_e + N_v) \}.
\]
%Recall that $N_e$ is the number of coarse-grid edges and $N_v$ is the number of coarse-grid nodes.
One can reduce the cost of solving local problems by using the randomization techniques
\cite{randomized2014}.

\subsubsection{Offline Space}

\def\BBPhi{{\Phi}}
\def\BBPsi{{ \Psi}}

We  perform a space reduction in the snapshot space through
using a local spectral problem in $\omega_i^{\eps}$.
We consider the following local eigenvalue problem in the snapshot space
\begin{equation}\label{eq:offeig2}
A^{\omega_i^\epsilon} \BBPsi_k^{\omega_i^\epsilon} = \lambda_k^{\omega_i^\epsilon} S^{\omega_i^\epsilon}\BBPsi_k^{\omega_i^\epsilon},
\end{equation}
where
$ (A^{\omega_i^\epsilon})_{mn}
= a_i(\psi_m^{\omega_i^{\eps}} , \psi_n^{\omega_i^{\eps}} )$,
$(S^{\omega_i^{\eps}})_{mn}
= s_i(\psi_m^{\omega_i^{\eps}} , \psi_n^{\omega_i^{\eps}} )$
and
$a_i(u, v) = \int_{\omega_i^{\eps}} \nabla u : \nabla v $,
$s_i(u, v) = \int_{\omega_i^{\eps}} |\nabla \chi_i |^2 u \cdot v$
with  $\chi_i$  specified below.
Note that we solve the above spectral problem in the local snapshot space
corresponding to the neighborhood domain $\omega_i^{\eps}$.
We arrange the eigenvalues in an increasing order,
choosing the first $l_i$ eigenvalues and taking the corresponding eigenvectors
$\BBPsi_k^{\omega_i^\epsilon}$, for $k = 1,2,...,l_i$,
to form the basis functions,
i.e.,
$\widetilde{\BBPhi}_k^{\omega_i^\epsilon} = \sum_j \Psi_{kj}^{\omega_i^\epsilon} \psi_j^{\omega_i^\epsilon}$, where $\Psi_{kj}^{\omega_i^\epsilon}$ are the coordinates of the vector $\BBPsi_{k}^{\omega_i^\epsilon}$.
Then we define
\[
\widetilde{V}^{\omega_i^{\eps}}_{H,\text{off}} = \text{span} \{ \widetilde{\BBPhi}_k^{\omega_i^\epsilon}, ~~k = 1,2,..., l_i \}. \numberthis\label{locofft}
\]

For constructing the conforming offline space,
we multiply the functions $\widetilde{\BBPhi}_k^{\omega_i^\epsilon}=(\widetilde{\Phi}_{x_1, k}^{\omega_i^\epsilon},\widetilde{\Phi}_{x_2, k}^{\omega_i^\epsilon})$
by a partition of unity function $\chi_i$.
Note that discontinuous Galerkin methods will eliminate the need for
$\chi_i$.
We remark that we define the partition of unity functions $\{ \chi_i\}$ with respect to the coarse nodes
and the mid-points of coarse edges. One can choose $\{ \chi_i\}$ to be the standard multiscale finite element basis.
However, upon multiplying by partition of unity functions, the resulting basis functions
no longer  have constant divergence,
which affects the scheme's stability.
Moreover, the multiplication by partition of unity functions may not honor
the perforation boundary conditions as discussed before.
To resolve the problem with the divergence, for each $\widetilde{\BBPhi}_k^{\omega_i^\epsilon}$ we
 find two functions  $  \BBPhi_{1, k}^{\omega_i^\epsilon}$ and $ \BBPhi_{2, k}^{\omega_i^\epsilon}$ (these are vector functions) that solve two local optimization
problems in every coarse-grid element $K^i_j \subset \omega_i^{\eps}$:
\begin{equation}\label{hrm-x}
\min{ \norm{\nabla \BBPhi_{1, k}^{\omega_i^\epsilon}}_{L^2(K^i_j)} } \text{~~such that~~}
\text{div}( \BBPhi_{1, k}^{\omega_i^\epsilon})
=  \frac{1}{|K^i_j|} \int_{\partial K^i_j }  (\chi_i \widetilde{\Phi}_{x_1, k}^{\omega_i^\epsilon} , 0) \cdot n_{ij} \, ds \quad \text{in } K^i_j
\end{equation}
with $n_{ij}$ the outward normal vector to $ \partial K^i_j$ and the boundary condition
$\BBPhi_{1, k}^{\omega_i^\epsilon}  =   (\chi_i \widetilde{\Phi}_{x_1, k}^{\omega_i^\epsilon} , 0), \text{ on } \partial K^i_j$;
and
\begin{equation}\label{hrm-y}
\min{ \norm{\nabla \BBPhi_{2, k}^{\omega_i^\epsilon}}_{L^2(K^i_j)} } \text{~~such that~~}
\text{div}( \BBPhi_{2, k}^{\omega_i^\epsilon})  =
\frac{1}{|K^i_j|} \int_{\partial K^i_j }  (0, \chi_i \widetilde{\Phi}_{x_2, k}^{\omega_i^\epsilon} ) \cdot n_{ij} \, ds \quad \text{in } K^i_j,
\end{equation}
with the boundary condition
$\BBPhi_{2, k}^{\omega_i^\epsilon}  =    (0, \chi_i \widetilde{\Phi}_{x_2, k}^{\omega_i^\epsilon} ), \text{ on } \partial K^i_j$.

Combining them, we obtain the global offline space:
\[
V_{H,\text{off}}  = \text{span} \{
\BBPhi_{1, k}^{\omega_i^{\epsilon}} \text{ and } \BBPhi_{2, k}^{\omega_i^\epsilon}: \quad
1 \leq i \leq (N_e+N_v)  \text{  and  }   1 \leq k \leq l_i
\}.
\]
Using a single index notation, we can write
$ V_{H,\text{off}} = \text{span} \{ \BBPhi^{\text{off}}_{i} \}_{i=1}^{N_u}$,
where $N_u =\sum_{i=1}^{N_e+N_v} l_i$.
We will use this  as an approximation space for the velocity.
For coarse approximation of the pressure, we will take $Q_{H,\text{off}}$
to be the space of piecewise constant functions on the coarse mesh.
Note that for many applications,
the size of the coarse problem can be substantially small compared
to the fine-scale problem. For example, for problems with scale separation,
we can only use 1-2 basis functions.
Once multiscale basis functions are constructed,
the global coarse-grid problem can be solved on the coarse grid with pre-computed
multiscale basis functions.

\subsubsection{Numerical results}
\label{sec:numerical}

\begin{figure}
\centering
\includegraphics[width=5.5in, height=2.5in]{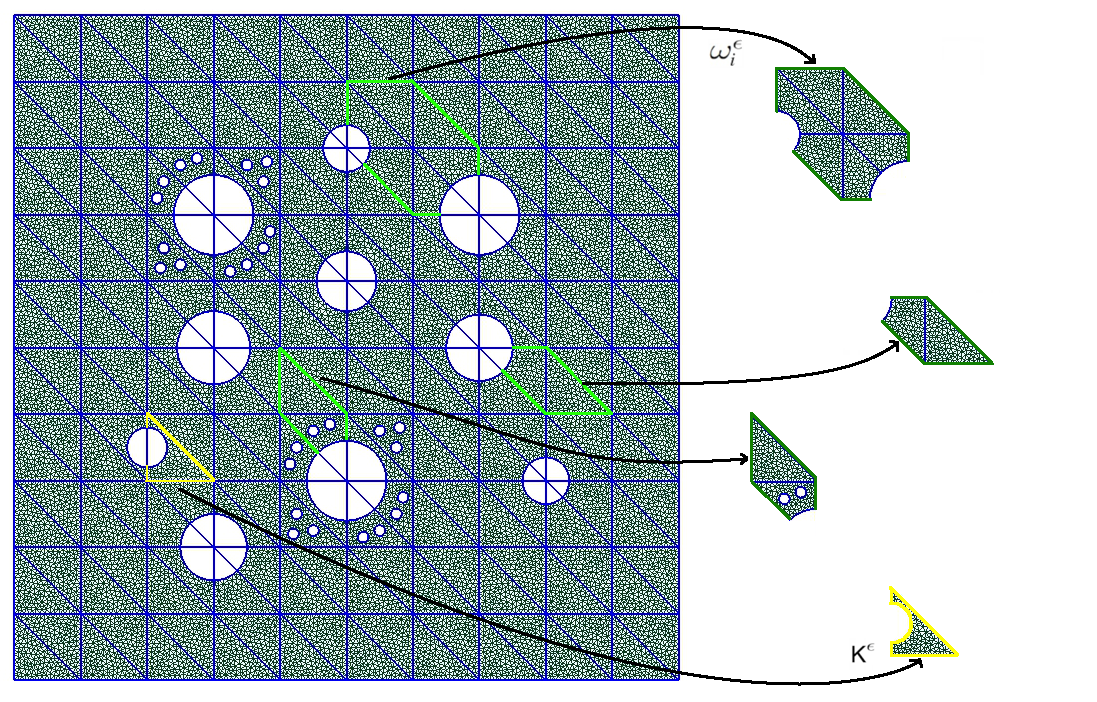}
\caption{Heterogeneous perforated medium  and the description of coarse and fine meshes.}
\label{fig:twodomain}
\end{figure}

We present a numerical result
for Stokes equations in $\Omega = [0,1] \times [0,1]$.
We consider a perforated domain as illustrated in Figure
\ref{fig:twodomain}, where
the perforated regions $\mathcal {B}^{\eps}$ are circular.
We coarsely discretize the computational domain  using uniform triangulation, where the
coarse mesh size $H=\frac{1}{5}$ for Stokes problem.
Note that our  approach does not impose any geometry restriction on coarse meshes.
 Furthermore, we use nonuniform triangulation inside each coarse triangular
element to obtain a finer discretization.

\begin{figure}[!h]
\begin{center}
        \includegraphics[width=1.0\textwidth]{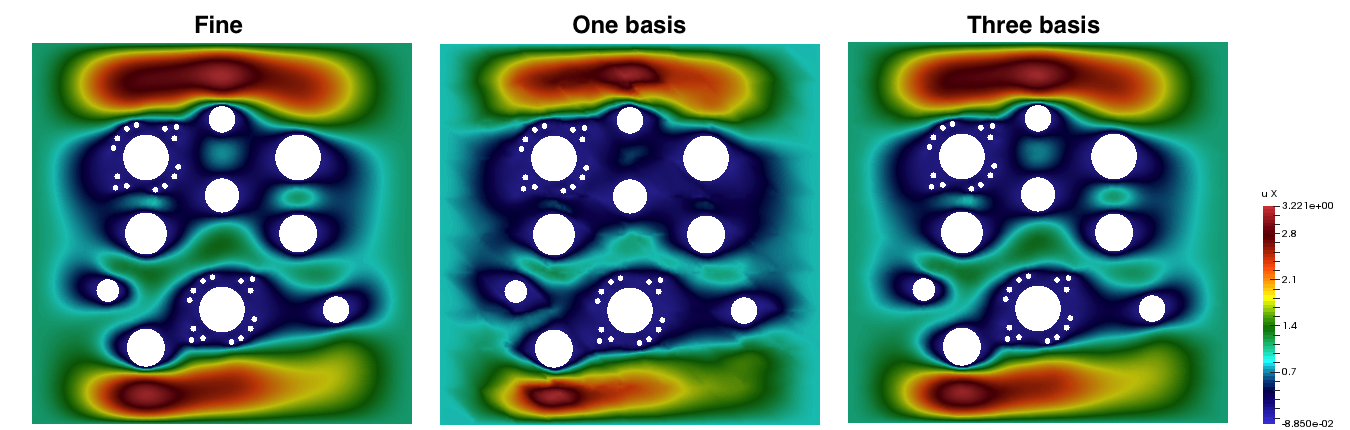}
\end{center}
    \caption{Stokes problem. Fine-scale and multiscale solutions for velocity
    $u_1$ in perforated domain (Figure \ref{fig:twodomain}).
    Left: fine-scale solution, $DOF = 171688$.
    Middle: multiscale solutions using 1 multiscale basis function for the velocity ($L_2$ error is $10.1$\%).
    Right: multiscale solutions with 3 multiscale basis function for the velocity ($L_2$ error is $1.02$\%). }
    \label{fig:st-ex}
\end{figure}

We consider the Stokes operator
with zero velocity $u = (0, 0)$ on $\partial \om^{\eps}\cap  \partial \mathcal {B}^{\eps}$
and ${\partial u\over\partial n} = (0, 0)$ on $\partial \om$.
For the fine-scale approximation of the Stokes problem, we use
$P^2$-elements for velocity and piecewise constants for pressure.
The $x_1$-component of the velocity solution is shown in Figure \ref{fig:st-ex}.
We present the offline solution with one
basis function and three basis functions.
We observe that the offline velocity solution
with one basis function (DOF=1082)
is not able to capture the solution accurately. The
$L^2$-error is $10.1$\%.  If we use three multiscale basis functions
per coarse element (DOF=2846),
we obtain a more accurate solution,  the $L^2$-error is only $1.02$\%.
These results show that one needs a systematic and rigorous approach for enriching
coarse-grid spaces in order to obtain  accurate solutions using only a small number of  basis
functions.

\section{Selected Applications}
\label{sec:applications}

Previously, we discussed the main concept of the GMsFEM and some important
ingredients. The GMsFEM can be used in various applications.
Below, we discuss a few applications.

\subsection{Two-phase flow}

We present simulation results for two-phase flow and
transport problems. We consider the two-phase
 flow problem with zero Neumann
boundary condition
\begin{equation}
\label{eq:twophase_flow}
\begin{split}
-\eta(S)\kappa\nabla p = {v} \mbox{ in }\Omega,\
\text{div} \, {v} = f \mbox{ in }\Omega,\
v\cdot n = 0 \mbox{ on }\partial \Omega,
\end{split}
\end{equation}
where
$\eta(S) = \frac{\kappa_{rw}(S)}{\mu_w}+\frac{\kappa_{ro}(S)}{\mu_o}$
and
$\kappa_{rw}(S) = S^2$,
$\kappa_{ro}(S)= (1-S)^2$, $\mu_w = 1$,
$\mu_o = 5$.
The saturation equation is given by
\begin{equation*}
S_t + v \cdot \nabla F(S) = r,
\end{equation*}
where
$F(S) =
\frac{\kappa_{rw}(S)/\mu_w}{\kappa_{rw}(S)/\mu_w+\kappa_{ro}(S)/\mu_o}$.
The above flow equation (\ref{eq:twophase_flow})
is solved by the mixed GMsFEM, and the
saturation equation is solved on the fine grid by the finite volume
method.
In our simulations, we take $f$ to be zero except for the
top-left and bottom-right fine-grid elements, where $f$ takes the
values of $1$ and $-1$, respectively. Moreover, we set the initial
value of $S$ to be zero. For the source $r$, we also take it as zero
except for the top-left fine element where $r=1$.

\begin{figure}[htp!]
\centering \subfigure[Fine-grid solution]{
\includegraphics[width=2in, height=1.5in]{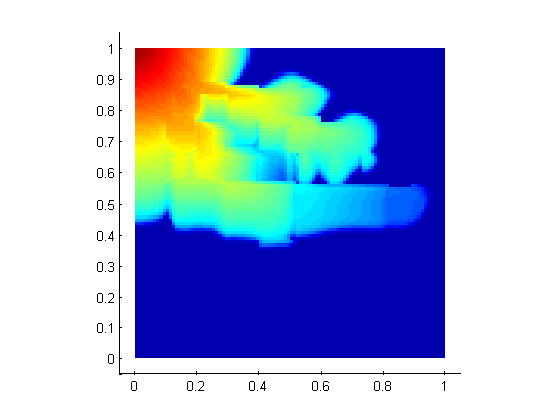}}
\subfigure[Relative $L^2$ error = 1.8\%]{
\includegraphics[width=2in, height=1.5in]{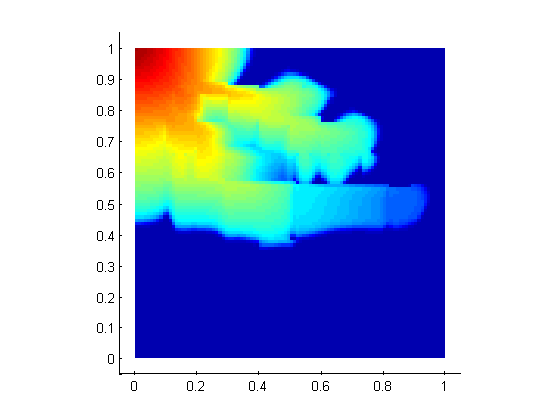}}
%\subfigure[Relative $L^2$ error = 1.8\%]{
%\includegraphics[scale=.23]{2off3_t5000.eps}}
%\includegraphics[width=1.5in, height=1.in]{2off3_t5000.eps}}
%\includegraphics[width=1.5in, height=1.in]{S_off_3000.eps}}
  \caption{Saturation solution obtained by using $v_o$ ($10\times10$ coarse grid, 3 basis per coarse edge)}\label{fig:10off4t}
\end{figure}

For the simulations of two-phase flow, the mixed GMsFEM
is used for the flow equations.
The multiscale basis functions are computed at time zero
using the unit mobility, i.e., $\eta=1$. These multiscale basis function
for the velocity field is used without any modification
to compute the fine-scale velocity field. The fine-scale
velocity field is further used to update the saturation file.
We note that each time the flow equation is solved on a coarse grid, which
provides a substantial computational saving.
We present the numerical results for the saturation in
Figure \ref{fig:10off4t} for the permeability field shown in Figure
\ref{fig:perm_random}. As we observe from this figure that with three
basis functions, we obtain very good agreement ($L^2$ error is $1.8$\%).

\subsection{Flow in fractured media}

\begin{figure}[htb]
  \centering
 \includegraphics[width=6in, height=2in]{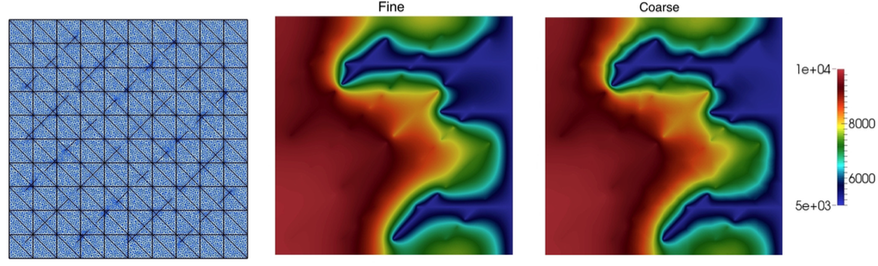}
  \caption{Left: the fracture distributions. Middle: the fine-scale solution.
Right: the coarse-grid solution.}
  \label{fig:shale_results}
\end{figure}

In this section, we briefly discuss the application of the GMsFEM to flows
in fractured media, in particular, the application to shale gas transport
in fractured media.
We remark that the study of flows in fractured media has a long
history (see, e.g.,  \cite{Lee01, hkj12}).
We study nonlinear
gas transport in fractured media
(cf., \cite{aes14})
\begin{equation}
\label{eq:model1}
a_{z}(c) {\partial c \over \partial t} = \text{div}(b_{z}(c,x) \nabla c),
\end{equation}
where $z=m,f$ with $m$ referring to the matrix phase and $f$ referring to
the fracture phase. The fractures have ``zero'' thickness and are
represented on a fine-scale variational formulation as edge elements.
More precisely, the domain $\Omega$ can be represented by
$\Omega = \Omega_m \oplus_i \, d_i \Omega_{f,i}$.
Here, $d_i$ is the aperture of the $i$ th fracture and $i$ is the index of the fractures (see Figure \ref{fig:shale_results}). Then, the fine-scale equations have the form
\begin{equation}
\label{eq:discr1}
\begin{split}
%m( \frac{\partial c}{\partial t}, v) +
%a(c, v) &=
\int_{\Omega_m}   a_{m}(c) {\partial c \over \partial t} v  +
 \sum\limits_i
d_i\int_{\Omega_{f,i}} a_{f}(c) {\partial c \over \partial t} v +
+\int_{\Omega_m}   b_{m}(c,x)\nabla c \cdot \nabla v  +
 \sum\limits_i
d_i\int_{\Omega_{ f,i}} b_{m}(c,x)\nabla c\cdot \nabla v  = 0.
\end{split}
\end{equation}

In our numerical simulations, we solve the above problem using the GMsFEM
framework for some specific parameter values that are specified below.
We note that the fractures are modeled within snapshot space
via harmonic extensions as we discussed above (see \cite{akkutlu2015multiscale}
for more details). In the method, the basis functions
are constructed using the steady state with $b_z=1$ in (\ref{eq:model1})
and taking into account the fracture distribution. The
problem (\ref{eq:model1}) is solved by using implicit discretization
and by linearizing $c$ at the previous time step.

%We consider the solution of problem with constant  and
%nonlinear matrix-fracture coefficients in \eqref{eq:discr1}. As constant coefficients
%(see previous section) representing matrix and fracture properties,
We use following parameters:
$a_{{m}}(c) 		= \phi 	+ (1-\phi) {\partial F}/{\partial c}$;
$b_{m}(c,x) 	= \phi D + (1-\phi) D_s {\partial F}/{\partial c} + \phi {\kappa} R T c/\mu$; $a_{f}(c)  = \phi_f$;
$b_{f}(c,x) = {\kappa_f} R T c/\mu$;
$D = 5 10^{-7} [m^2/s]$, $\phi = 0.04$,
$T = 413 [K]$, $\mu = 2 \cdot 10^{-5} [kg/(m \, s)]$ and for fractures $k_f = 10^{-12} [m^2]$, $\phi_f = 0.001$.
As for permeability $\kappa$, we use constant $\kappa= 10^{-18} [m^2]$.
$p = RT c$, $p_c = 10^9 [Pa]$, $p_1 = 1.8 \cdot 10^9 [Pa]$, $\alpha = 0.5$ and
$\mathcal{M} = 0.5$.
For the sorbed gas, we use Langmuir model
$F(c)= c_{\mu s} \frac{s}{(1 + s c)^2}$,
where $s = 0.26 \cdot 10^{-3}$ and $c_{\mu s} = 0.25 \cdot 10^{-5}  [mol/m^3]$.
%$a_m = 0.8$, $b_m = 1.3 \cdot 10^{-7}$,
%$a_f = 0.001$, $b_f = 1.0$.
%For nonlinear matrix-fracture coefficients, we use
%$a_{f}(c)  = \phi_f$,
%$b_{f}(c,x) = \frac{\kappa_f}{\mu} R T c$,
%$D_k = 10^{-7} [m^2/s]$, $D_i = 10^{-8} [m^2/s]$, $\phi = 0.04$, $T = 413 [K]$, $\mu = 2 \cdo%t 10^{-5} [kg/(m \, s)]$ and for fractures $k_f = 10^{-12} [m^2]$, $\phi_f = 0.001$.
We present the numerical results in Figure \ref{fig:shale_results}.
In the left figure, we depict the fracture distributions, and
then, show the fine-scale and the coarse-scale solution.
We observe that with only $4$ basis functions, we obtain a good
agreement between the fine-scale and multiscale solution ($L^2$ error $0.34$\%,
$H^1$ error $4.69$\%).

\subsection{Other applications}

We have demonstrated two applications of the GMsFEM.
The GMsFEM can be used for other complex multiscale problems.
For example, in \cite{chung2015_elasticwave_fractured},
we have applied the GMsFEM for elastic
wave propagation in fractured media. We have also used the GMsFEM
for acoustic wave propagation \cite{eric-2012}. The GMsFEM can also be effectively used
in inverse problem, where the forward problem is solved many times.
In \cite{efendiev2015multilevel},
we apply the GMsFEM in uncertainty quantification in
inverse problems.

\section{Discussions}
\label{sec:summary}

In this paper, we give an overview of multiscale finite element methods
and discuss their generalizations. Due to page limitations,
we do not present a comprehensive
 overview of multiscale methods and mention
only some important aspects of the GMsFEM as a tool for performing adaptive
multiscale model reduction.
%We present a general methodology
%for performing adaptive multiscale model reduction and some
%important ingredients for the method.
The application of the method to
nonlinear problems is not discussed in the paper. We
refer to
(\cite{chung2015_nonlinear}) for the application
to nonlinear problems, where we consider
nonlinear eigenvalue problems, hybridization, and nonlinear interpolation
ideas to efficiently solve these problems.

The proposed local reduced-order models can be used in performing global model
reduction. The GMsFEM can be employed to construct
the snapshots for the global model reduction methods. In this case,
using reduced representation of the global snapshots, we can perform
POD or another model reduction for the global modes. Moreover, we
can use adaptively local approaches to update the global
snapshots when online computations are needed. These global-local
approaches can be effective in many applications \cite{alotaibi2015global}.

Optimal multiscale basis functions can be obtained using the oversampling techniques, optimization,
and singular value decomposition
\cite{hou2015optimal}. In this paper, we focus on a general adaptive strategy that can be
easily adapted for
various applications and discretization, as demonstrated.
One can use the ideas of optimized multiscale
basis functions as in \cite{hou2015optimal} in the proposed framework.

The stochastic multiscale problems are challenging due to additional degrees of freedom
due to randomness.
For stochastic problems, one can use multiscale finite element framework
with Monte Carlo techniques or separation of variables. A promising approach
for solving multiscale stochastic problems is to use data driven stochastic
method concepts  as shown in \cite{zhang2015multiscale}.

Multiscale framework proposed in the paper can also be employed for the stabilization.
For example, in \cite{calo2015multiscale}, the authors use the GMsFEM to stabilize
the convection-dominated problems by constructing multiscale test spaces.

Finally, we would like to remark that the proposed methods are well suited for
the parallel computation,
in particular, the task of multiscale basis computations is embarrassingly parallel.
Our numerical results (using up to $1000$ processors) on 3D perforated problems
show that one can achieve almost perfect scaling in
the parallel computations.

{\bf Acknowledgements.} TYH's research was in part supported by DOE Grant DE-FG02-06ER257, and NSF Grant No. DMS-1318377, DMS-1159138.
EC's research was in part supported by Hong Kong RGC General Research Fund Project No. 400813.

\bibliographystyle{plain}
\bibliography{references,references1}
%\bibliography{references,references1,Fractures,rlg,seismics,SeismicWaves,SeismicAnisotropy,SeismicProcessing,Inversino,refs,chester,YE_2015_references,mlmc}

\end{document}